\documentclass[12pt]{article}
\usepackage{mathrsfs}
\usepackage{amsfonts}
\usepackage{amsmath, amsthm, amsfonts, amssymb}
\theoremstyle{definition}

\numberwithin{equation}{section} \theoremstyle{remark}

\def\<{\langle}
\def\>{\rangle}

\def\ra{\rightarrow}

\def\a{\alpha}

\def\CC{{\bf C}}
\def\RR{{\bf R}}
\def\ZZ{{\bf Z}}

\def\-{\overline}

\def\endpf{\hbox{\vrule height1.5ex width.5em}}

\def\b{\beta}
\def\a{\alpha}
\def\endpf{\hbox{\vrule height1.5ex width.5em}}
\def\RR{{\Bbb R}}
\def\CC{{\Bbb C}}

\def\ZZ{{\Bbb Z}}

\def\M*{\wt{M^*}}

\def\-{\overline}

\def\wt{\widetilde}
\def\ra{\rightarrow}

\def\endpf{\hbox{\vrule height1.5ex width.5em}}

\def\a{\alpha}

\def\b{\beta}

\def\a{\alpha}

\def\endpf{\hbox{\vrule height1.5ex width.5em}}
\begin{document}
\bigskip
\title{\bf A codimension two CR singular submanifold  that is formally equivalent to a symmetric quadric}
\author{Xiaojun Huang\footnote{
Supported in part by NSF-0500626} \ \  and\  Wanke Yin}
\date{ }

\def\beq{\begin{equation}}
\def\nneq{\end{equation}}
\def\beqn{\begin{eqnarray}}
\def\neqn{\end{eqnarray}}
\def\beqna{\begin{eqnarray*}}
\def\neqna{\end{eqnarray*}}
\def\bedis{\begin{displaymath}}
\def\nedis{\end{displaymath}}
\def\-{\overline}

\vspace{6cm} \maketitle

\begin{abstract}
Let  $M\subset \mathbb{C}^{n+1}$ ($n\geq 2$) be a real
 analytic submanifold defined by an equation of the form: $w=|z|^2+O(|z|^3)$, where we use $(z,w)\in \mathbb{C}^{n}\times \mathbb{C}$
for the coordinates of $\mathbb{C}^{n+1}$.
 We first derive a pseudo-normal form for $M$ near $0$. We then use it to prove
  that $(M,0)$ is holomorphically equivalent to the quadric $(M_\infty: w=|z|^2,\ 0)$
 if and only if it can be formally transformed to $(M_\infty,0)$. We also use it to give a necessary and sufficient condition
when
 $(M,0)$ can be formally flattened.
 The result is due to  Moser for the case of  $n=1$.
\end{abstract}

\section{Introduction}
Let $M\subset \mathbb{C}^{n+1}$ ($n\ge 1$) be a submanifold. For a
point $p\in M$, we define $CR(p)$ to be the CR dimension of $M$ at
$p$, namely, the complex dimension of the space $T^{(0,1)}_pM$. A
point $p\in M$ is called a CR point if $CR(q)=CR(p)$ for $q(\in
M)\approx p$. Otherwise, $p$ is called a CR singular point of $M$.
The local equivalence problem in Several Complex Variables is to
find a complete set of holomorphic invariants of $M$ near a fixed
point $p\in M$. The investigation normally has  quite different
nature  in terms of  whether $p$ being a CR point or a CR singular
point. The  CR case was first considered by Poincar\'e and Cartan. A
complete set of invariants in the
 strongly pseudoconvex hypersurface  case was give by Chern-Moser in [CM] (see
the survey paper of Baouendi-Ebenfelt-Rothschild [BER1] and the
lecture notes of the first author [Hu1] for many references along
these lines).
The study for the CR singular points first appeared   in the paper
of Bishop [Bis]. Further investigations on the precise holomorphic
structure of $M$ near a non-degenerate CR singular point, in the
critical dimensional case of $dim_{\RR}M=n+1$, can be found in the
work of Moser-Webster [MW] and in the work of   Moser [Mos], Gong
[Gon1-2], Huang-Yin [HY], etc. (The reader can find many references
in [Hu1] on this matter.)


Recently, there  appeared several papers,  in which CR singular
points in the non-critical dimensional case were considered (see
[Sto], [DTZ], [Cof1-2], to name a few). In [Sto], among other things,
Stolovitch introduced a set of generalized Bishop invariants for a
non-degnerate general CR singular point, and established some of the
results of Moser-Webster [MW]  to the case of
$dim_{\RR}M>dim_{\CC}{\CC}^{n+1}$. In [DTZ],
Dolbeault-Tomassini-Zaitsev introduced the concept of the elliptic
flat CR singular points and studied global filling property by
complex analytic varieties for a class of compact submanifold  of
real codimension two in ${\CC}^{n+1}$ with exactly two  elliptic
flat CR singular points.

In this paper,  we study the local holomorphic structure of  a
manifold $M$ near a CR singular point $p$, for which
 we can  find a local holomorphic change of coordinates
such that in the new coordinates system,  $p=0$ and $M$ near $p$ is
defined by an equation of the form: $w=|z|^2+O(|z|^3)$. Here we use
$(z,w)\in {\CC}^{n}\times \CC$ for the coordinates of ${\CC}^{n+1}$.
Such a non-degenerate CR singular point  has an intriguing nature
that its quadric model  has the largest possible symmetry. We will
first derive a pseudo-normal form for $M$ near $p$ (see Theorem
2.3). As expected,  the holomorphic structure of $M$ near $p$ is
 influenced  not only by
the nature of the CR singularity, but also by the fact that $(M,p)$
partially inherits the property of strongly pseudoconvex CR
structures for $n>1$.
 Unfortunately, as in the case of
$n=1$ first considered by Moser [Mos], our pseudo-normal form is
still subject to the simplification of the complicated infinite
dimensional formal automorphism group of the quadric
$aut_0(M_\infty)$, where $M_\infty$ is defined by $w=|z|^2$. Thus,
our pseudo-normal form can not be used to solve the local
equivalence problem. However, with the rapid iteration procedure, we
will show in $\S 4$ that if all higher order terms in our
pseudonormal form vanish,
 then $M$ is biholomorphically equivalent
to the model $M_\infty$. Namely, we have the following:

\bigskip
{\bf Theorem 1}: {\it Let  $M\subset \mathbb{C}^{n+1}$ ($n\ge 1$) be
a real  analytic submanifold defined by an equation of the form:
$w=|z|^2+O(|z|^3)$. Then $(M,0)$ is holomorphically equivalent to
the quadric $(M_\infty, 0)$ if and only if it can be formally
transformed to $(M_\infty,0)$.}
\bigskip

One of the  differences of our consideration here from the
case of $n=1$ is that a generic $(M,0)$ can not be formally mapped into
the Levi-flat hypersurface $Im(w)=0$. As another  application of
the pseudo-normal form to be obtained in $\S 2$, we will  give
a necessary and sufficient condition when $(M,0)$ can be
formally  flattened (see Theorem 3.5).

Theorem 1, in the case of $n=1$, is due to Moser [Mos]. Indeed, our
proof of Theorem 1  uses  the  approach of  Moser in [Mos] and Gong
in [Gon2], which is   based on the rapidly convergent power series
method. Convergence results along the lines of Theorem 1 near other
type of   CR singular points
 can be found in the earlier papers of
Gong [Gon1] and Stolovitch [Sto]. The papers of Coffman [Cof1-2] also contain the
rapid convergence arguments in the setting of other CR singular cases.

\section{A formal pseudo-normal form}

We use $(z,w)=(z_1,\cdots,z_n,w)$ for the coordinates in
$\mathbb{C}^{n+1}$ with $n\ge 2$ in all that  follows. We first
recall some notation and definitions already discussed  in the
previous papers of Stolovitch [Sto] and Dolbeault-Tomassini-Zaitsev
[DTZ].

 Let
$(M,0)$ be a formal submanifold of codimenion two in ${\Bbb
C}^{n+1}$ with $0\in M$ as a CR singular point and
$T^{(1,0)}_0M=\{w=0\}$. Then, $M$ can be defined by a formal
equation of the  form:
\begin{equation}
w=q(z,\-{z})+o(|z|^2), \label{eqn:010}
\end{equation}
 where $q(z,\-{z})$ is a quadratic polynomial in
$(z,\-{z})$. We say that $0\in M$ is a not-completely-degenerate CR
singular point if there is no change of coordinates in which we can
make $q\equiv 0.$ Following Dolbeault-Tomassini-Zaitsev, we further
say that $0$ is a not-completely-degenerate flat CR singular point
if we can make $q$ real-valued after a  linear change of
variables.

Assume that $0$ is a not-completely-degenerate flat CR singular
point with $q(z,z)=A(z,\-{z})+B(z,\-{z})\in {\Bbb R}$  for each $z$.
Here $A(z,\-{z})=\sum_{\a,\b=1}^na_{\a\-{\b}}z_\a\-{z_\b},\
B(z,\-{z})=2Re(\sum_{\a,\b=1}^nb_{\a\b}z_az_\b).$ Then the
assumption that $A(z,\-{z})$ is definite is independent of the
choice of the coordinates system. Suppose that $A$ is definite. Then
making use of the classical Takagi theorem, one can find a linear
change of coordinates in $(z,w)$ such that in the new coordinates,
in the defining equation for $(M,0)$ of the form in (\ref{eqn:010}),
one has that
$$q(z,\-{z})=\sum_{\a=1}^n\{|z_\a|^2+\lambda_\a
(z_\a^2+\-{z_\a}^2)\},$$ where $0\le \lambda_\a<\infty$ with
$0\le\lambda_1\le\cdots \le \lambda_n< \infty$. In terms of
Stolovitch, we call $\{\lambda_1,\cdots,\lambda_n\}$ the set of generalized
Bishop invariants. When $0\le \lambda_\a<1/2$ for all $\a$,
 we say that $0$ is an
elliptic flat CR singular point of $M$. Notice that $0\in M$  is an
elliptic flat CR singular point if and only if in a certain defining
equation of $M$ of the form as in (\ref{eqn:010}), we can make
$q(z,\-{z})>0$ for $z\not = 0$. (Hence the  definition  coincides
with the notion of elliptic flat Complex points  in [DTZ].)
 When
$\lambda_\a>1/2$ for all $\a$, we say $0\in M$ is a hyperbolic flat
CR singular point. Notice that, in the other case, we can always
find a two dimensional linear subspace of ${\Bbb C}^{n+1}$ whose
intersection with $M$ has a parabolic complex tangent at $0$. For a
more general related notion on  ellipticity and hyperbolicity, we
refer the reader to the paper of Stolovitch [Sto].

In terms of the terminology above,  the manifold  in Theorem 1
 has  vanishing
generalized Bishop invariants at the CR singular  point. In [Gon1]
 [Sto],
 one finds
 the study on the related  convergence problem in the other
situations, where, among other non-degeneracy conditions,
 all the generalized Bishop invariants are assumed to be
non-zero. The
method studying CR singular points with
vanishing Bishop invariants is
different from that used in the non-vanishing Bishop invariants case
(see [MW] [Mos] [Gon2] [Sto] [HY]).

We now return to the manifolds with only vanishing
generalized Bishop invariants.

\medskip

Let $E(z,\bar{z})$ $\left(\mbox{respectively}, \ f(z,w)\right)$ be a
formal power series in $(z,\bar{z})\left(\mbox{respectively, in}\
(z,w)\right)$ without constant term. We say
$Ord\left(E(z,\bar{z})\right)\geq k$ if $E(tz,t\bar{z})=O(t^k)$.
Similarly, we say $Ord_{wt}\left(f(z,w)\right)\geq k$ if
$f(tz,t^2w)=O(t^k)$. Set the weight of $z,\bar{z}$ to be 1 and that
of $w$ to be 2. For a polynomial $h(z,w)$, we define its weighted
degree, denoted by $deg_{wt}h$, to be the degree counted in terms
the weighted system just given. Write $E^{(t)}(z,\bar{z})$ and
$f^{(t)}(z,w)$ for the sum of monomials with weighted degree $t$ in
the  expansion of $E$ and $f$ at $0$, respectively.

Write $u_k=\sum_{i=1}^{k}|z_i|^2$ for $1\le k\le n$ and
$v_k=\sum_{i=1}^{k-1}|z_i|^2-|z_k|^2$ for $2 \leq k \leq n$. We also
write $u=u_n=|z|^2$. In what follows, we make a convention that the
sum $\sum_{p=j}^{l}a_p$ is defined to be $0$ if $j>l$.

 We start with
the following elementary algebraic lemma:

\medskip
{\bf Lemma 2.1:}
$Span_{\CC}\{|z_1|^2,\cdots,|z_n|^2\}=Span\{u,v_2,\cdots,v_n\}$.
Moreover, for each index $i$ with $1\le i\le n$, $|z_i|^2$ can be
uniquely expressed as the following linear combination of $u,\
v_2,\cdots,\ v_n$:

\begin{equation}
\left\{
\begin{array}{l}
|z_1|^2=2^{1-n}\left(u+\sum\limits_{h=2}^{n}2^{n-h}v_h\right),\\
|z_i|^2=2^{-(n+1-i)}\left(u+\sum\limits_{h=i+1}^{n}2^{n-h}v_h-2^{n-i}v_i\right)\
\mbox{for}\  2\leq i\leq n.
\end{array}
\right. \label{2.3}
\end{equation}

\medskip
{\it Proof of Lemma 2.1:} By a direct computation, we have
\begin{equation*}
\begin{array}{l}
2^{1-n}\left(u+\sum\limits_{h=2}^{n}2^{n-h}v_h\right)
   =2^{1-n}\left(\sum\limits_{i=1}^{n}|z_i|^2
   +\sum\limits_{h=2}^{n}2^{n-h}(\sum\limits_{i=1}^{h-1}|z_i|^2-|z_h|^2)\right)\\
   \hskip 1cm =2^{1-n}\left((1+\sum\limits_{h=2}^{n}2^{n-h})|z_1|^2+
   \sum\limits_{j=2}^{n-1}(1+\sum\limits_{h=j+1}^{n}2^{n-h}-2^{n-j})|z_j|^2\right)\\
   \hskip 1cm =2^{1-n}(2^{n-1}|z_1|^2)=|z_1|^2;\\
2^{-(n+1-i)}\left(u+\sum\limits_{h=i+1}^{n}2^{n-h}v_h-2^{n-i}v_i\right)\\
 \hskip 1cm
=2^{-(n+1-i)}\left(\sum\limits_{i=1}^{n}|z_i|^2+\sum\limits_{h=i+1}^{n}2^{n-h}
  (\sum\limits_{j=1}^{h-1}|z_j|^2-|z_h|^2)-2^{n-i}(\sum\limits_{j=1}^{i-1}|z_j|^2-|z_i|^2)\right)\\
 \hskip 1cm
=2^{-(n+1-i)}\left(\sum\limits_{j=1}^{i-1}(1+\sum\limits_{h=i+1}^{n}2^{n-h}-2^{n-i})|z_j|^2
  +(1+\sum\limits_{h=i+1}^{n}2^{n-h}+2^{n-i})|z_i|^2\right.\\
 \hskip 1.3cm
\left.+\sum\limits_{j=i+1}^{n}(1+\sum\limits_{h=j+1}^{n}2^{n-h}-2^{n-j})|z_j|^2\right)
=|z_i|^2,\ \hbox{for}\ i\ge 2.
\end{array}
\end{equation*}
 Hence, we see that
 $\hbox{span}_{\CC}\{|z_1|^2,\cdots,|z_n|^2\}=\hbox{span}_{\CC}\{u,v_2,\cdots,v_n\}.$
 The uniqueness assertion in the lemma now is obvious. $\endpf$

\medskip
For a formal (or holomorphic) transformation $f(z,w)$ of
$({\CC}^{n},0)$ to itself, we write
\begin{equation}
\left\{
\begin{array}{l}
f(z,w)=\left(f_1(z,w),\cdots,f_n(z,w)\right),\\
f_k(z,w)=\sum_{(i_1,\cdots,i_n)}f_{k,(I)}(w)z^I,\ I=(i_1,\cdots,i_n) \ \hbox{and}\ z^I=z_1^{i_1}
\cdots z_n^{i_n}.
\end{array}
\right. \label{2.4}
\end{equation}

Let $E(z,\bar{z})$ be a formal power series with $E(0)=0$. We next
prove the following:

\bigskip
{\bf Lemma 2.2:}  {\it $E(z,\bar{z})$ has the following
expansion:
\begin{equation}
E(z,\bar{z})=\sum_{\{i_k \cdot j_k=0,\ k=1\cdots,n\}}
E_{(I,J)}(u,v_2,\cdots,v_n)z^I{\overline{z}}^J
            =\sum_{\{i_k \cdot j_k=0,\ k=1,\cdots, n\}}
E_{(I,J)}^{(K)}z^I{\overline{z}}^Ju^{k_1}{v_2}^{k_2}\cdots
{v_n}^{k_n}. \label{2.5}
\end{equation}
Here and in what follows, we  write $I=(i_1,\cdots,i_n)$,
$J=(j_1,\cdots,j_n)$, $K=(k_1,\cdots,k_n)$, $z^I=z_1^{i_1}\cdots
z_n^{i_n}$ and ${\overline{z}}^J=\-z_1^{j_1}\cdots \-z_n^{j_n}$.
Moreover, the coefficients $E_{(I,J)}^{(K)}$ are uniquely determined
by $E$.}

 \medskip

{\it Proof of Lemma 2.2:} Since $\{|z_i|^2\}_{i=1}^{n}$ and
$\{u,v_2,\cdots,v_n\}$ are the unique linear combinations of  each
other by Lemma 2.1, one sees the existence of the expansion in
(\ref{2.5}). Also, to complete the proof of Lemma 2.3, it suffices
for us to prove the following statement:
$$ \sum\limits_{(I,J,K)\in
A(N,N^*)}E_{(I,J)}^{(K)}z^I{\overline{z}}^J|z_1|^{2k_1}\cdots
|z_n|^{2k_n}=0\ \mbox{if and only if }\ E_{(I,J)}^{(K)}\equiv 0.
$$
Here, we define $A(N,N^*)=\{(I,J,K)\in
\mathbb{Z}^n\times\mathbb{Z}^n\times\mathbb{Z}^n,\ i_l \cdot j_l=0,\
i_l,j_l,k_l\geq 0\ \mbox{for}\ 1 \leq l \leq n,\
\sum_{l=1}^n(i_l+k_l)=N,\ \sum_{l=1}^{n}(j_l+k_l)=N^*\}$. Let
$P=(p_1,\cdots,p_n)$ and $Q=(q_1,\cdots,q_n)$ with $p_1,
\cdots,p_n,q_1,\cdots, q_n$ non-negative integers be such that
$|P|=N, |Q|=N^*$. We define $A(N,N^*;P,Q)=\{(I,J,K)\in A(N,N^*):\
\ i_l \cdot j_l=0,\ i_l,j_l,k_l\geq 0, i_l+k_l=p_l,\ j_l+k_l=q_l,\
\mbox{for}\ 1 \leq l \leq n\}.$ Now, suppose that
$\sum\limits_{(I,J,K)\in
A(N,N^*)}E_{(I,J)}^{(K)}z^I{\overline{z}}^J|z_1|^{2k_1}\cdots
|z_n|^{2k_n}=0$. We then get
$$
\sum\limits_{(I,J,K)\in A(N, N^*;P,Q)}E_{(I,J)}^{(K)}\equiv 0, \
\hbox{for each} \  P, \ Q\ \hbox{with}\ |P|=N,\ |Q|=N^*.
$$
We next claim that  there is at most one element in $A(N, N^*;P,Q)$.
Indeed, $ (I,J,K)\in A(N, N^*;P,Q)$ if and only if $ i_l+k_l=p_l,\
j_l+k_l=q_l,\ i_l\cdot j_l=0$, for  $1\leq l \leq n.$ Now, if
$i_l=0$, then $k_l=p_l$. Since $j_l=q_l-p_l\geq 0$, thus this
happens only when $q_l\geq p_l$. If $j_l=0$, then $k_l=q_l$. Since
$i_l=p_l-q_l\geq 0$, we see that this can only happen when $p_l\geq
q_l$. Hence, we see that $i_l, j_l$ are uniquely determined by $p_l$
and $q_l$ when $p_l\not =q_l$. When $p_l=q_l$, it is easy to see
that $i_l=j_l=0, \ k_l=q_l=p_l$. We thus conclude the argument for
the claim. This  completes the proof of Lemma 2.3.
$\endpf$

\medskip

We now let   $M \subset\mathbb{C}^{n+1}$ be  a formal submanifold defined by:\\
\begin{equation}
w=|z|^2+E(z,\bar{z}) \label{2.1}
\end{equation}
where $E$ is a formal power series in $(z,\bar{z})$ with $Ord(E)
\geq 3$.
We will subject (\ref{2.1})  to  the following formal  power series transformation in $(z,w)$: \\
\begin{equation}
\left\{
\begin{array}{ll}
z'=F=z+f(z,w) \ &\ Ord_{wt}(f) \geq 2   \\
w'=G=w+g(z,w)      & \ Ord_{wt}(g) \geq 3.
\end{array}
\right.
\label{2.2}
\end{equation}

Write $e_j\in {\ZZ}^n$ for the vector whose component is $1$ at the
$j^{\hbox{th}}$-position and is  $0$ elsewhere.
We next give a formal
pseudo-normal form for $(M,0)$ in the following theorem:
\bigskip

{\bf Theorem 2.3:}\hspace{0.2cm} {\it There exits a unique formal
transformation of the form in (\ref{2.2}) with the  normalization
\begin{equation}
\begin{cases} f_{i,(0)}(u)=0, \ 1\le i\le n;\\
f_{i,(e_j)}(u)=0\ \mbox{for} \ 1\leq j<i\leq
n;\\
\ f_{1,(e_1)}(u)=0,\
\hbox{Im}\left(f_{i,(e_i)}(u)\right)=0\ \mbox{for} \ 2\leq i \leq n,
 \end{cases}
 \label{2.6}
\end{equation}
 that transforms M  to a formal submanifold defined in the following
 pseudo-normal form:
\begin{equation}
w'=|z'|^2+\varphi(z',\overline{z'}). \label{2.7}
\end{equation}
Here $\varphi=O(|z'|^3)$ and  in the following unique expansion of
$\varphi$,

\begin{equation}
\varphi =\sum_{i_l\cdot j_l=0,l=1,\cdots,n}
\varphi_{(I,J)}z^I{\overline{z}}^J=\sum_{i_k \cdot j_k=0,\
k=1,\cdots, n}
\varphi_{(I,J)}^{(K)}z^I{\overline{z}}^Ju^{k_1}{v_2}^{k_2}\cdots
{v_n}^{k_n}. \label{2.80}
\end{equation}
we have,  for any $k \geq 0,\ l\ge 1$, $\tau\geq 2$, the following
normalization condition:
\begin{equation}
\begin{cases}
\varphi_{(0,0)}^{(\tau e_1)}=0\ ;\\
 Re(\varphi_{(0,0)}^{(le_1+e_i)})=0,\
\mbox{for} \ 2 \leq
i \leq n\ ;\\
 \varphi_{(e_i,e_j)}^{(le_1)}=0, \  \mbox{for} \ i>j\ ;
 \\
\varphi_{(I,0)}^{(le_1)}=\varphi_{(0,I)}^{(le_1)}=\varphi_{(0,I)}^{(ke_1+e_j)}=0,
\ \mbox{for} \ |I|\geq 1;
\\ \varphi_{(I,e_h)}^{(ke_1)}=0,  \ \mbox{for}
\ h\geq 1,|I|\geq 2, i_h=0;
\\
\varphi^{(0)}_{(0,I)}=\-{\varphi^{(0)}_{(I,0)}}, |I|>2.
\end{cases}
\label{2.8}
\end{equation}
}

{\it Proof of Theorem 2.3}: We need to prove that the following
equation, with unknowns in $(f,g,\ \varphi)$, can be uniquely solved
under the normalization conditions in (\ref{2.6}) and (\ref{2.8}):
\begin{equation}
\begin{array}{l}
w+g(z,w)=\sum\limits_{i=1}^{n}\big(z_i+f_i(z,w)\big)\left(\bar{z_i}+\overline{f_i(z,w)}\right)
       +\varphi\left(z+f(z,w),\overline{z}+\overline{f(z,w)}\right).
\end{array}
\end{equation}
 Collecting terms of degree $t$
in the above equation, we obtain for each $t
\geq 3$ the following:
\begin{equation}
\begin{array}{ll}
E^{(t)}(z,\bar{z})+g^{(t)}(z,u)
=2Re\sum\limits_{i=1}^{n}\Big(\overline{z_i}f_i^{(t-1)}(z,u)\Big)
+\varphi^{(t)}(z,\bar{z})+I^{(t)}(z,\bar{z}), \label{hu-02}
\end{array}
\end{equation}
where $I^{(t)}(z,\bar{z})$ is a homogeneous polynomial of degree $t$
depending only on $g^{(\sigma)}$, $f^{(\sigma-1)}$,
$\varphi^{(\sigma)}$ for $\sigma < t$. Thus, by an induction
argument,  we need only to uniquely  solve the following  equation
under the above given normalization:
\begin{equation}
\Gamma(z,\bar{z})+g(z,u)=2Re\left(\sum\limits_{i=1}^{n}\left(\overline{z_i}f_i(z,u)\right)\right)+
\varphi(z,\bar{z}). \label{2.9}
\end{equation}
Indeed, if we can uniquely solve (\ref{2.9}), then, we can start with (\ref{hu-02}) with $t=3$ and $\Gamma=E^{(3)}$. We then get $(F^{(2)}, G^{(3)}).$
Now, we transform $M$ by $H_2=(z,w)+(F^{(2)}, G^{(3)})$. Then the new manifold is normalized up to weighted order $3$. Let $H=(F,G)=(z+O_{wt}(3),w+O_{wt}(4))$
be a normalized map and
 consider (\ref{hu-02}) with $t=4$. We can then uniquely determine $(F^{(3)},G^{(4)})$.
 Transforming the manifold by the map
  $H_2=(z,w)+(F^{(3)}, G^{(4)})$, we get  one which is normalized up to order $4$. Now, by an induction, we can prove the existence part of
 Theorem 2.3. The uniqueness part of the Theorem follows also from the unique solvability of (\ref{2.9}).

Expand $\Gamma$, $\varphi$ as in (\ref{2.5}) and (\ref{2.80}) and
expand $f$, $g$ as in  (\ref{2.4}). Making use of Lemma 2.2 and
comparing the coefficients in (\ref{2.9}) of $z^I\-{z}^J$ with
$i_l\cdot j_l=0, \ l=1,\cdots, n$, we get the following system:
\begin{align}
z^0\-{z}^0:&\ \ \ -g_{(0)}+\sum\limits_{i=1}^{n}2Re\left(|z_i|^2f_{i,(e_i)}\right)+\varphi_{(0,0)}=\Gamma_{(0,0)};\label{3.2.10}\\
z_j,\ \-{z_j}:&\ \ \ \begin{cases}
   -g_{(e_j)}+\overline{f_{j,(0)}}+\sum\limits_{i=1}^{n}|z_i|^2f_{i,(e_i+e_j)}+\varphi_{(e_j,0)}=\Gamma_{(e_j,0)}\\
   f_{j,(0)}+\sum\limits_{i=1}^{n}|z_i|^2\overline{f_{i,(e_i+e_j)}}+\varphi_{(0,e_j)}=\Gamma_{(0,e_j)}
      \end{cases}
            \hskip5pt \mbox{for}\ 1\leq j\leq n;\label{3.2.11}\\
z_i\-{z_j}:&\ \ \  \begin{cases}
   f_{j,(e_i)}+\overline{f_{i,(e_j)}}+\varphi_{(e_i,e_j)}=\Gamma_{(e_i,e_j)}\\
   \overline{f_{j,(e_i)}}+f_{i,(e_j)}+\varphi_{(e_j,e_i)}=\Gamma_{(e_j,e_i)}
   \end{cases} \hskip5pt \mbox{for}\ i\neq j;\label{3.2.12}\\
z_i\-{z}^J,\ z^J\-{z_i}:&\ \ \  \begin{cases}
    \overline{f_{i,(J)}}+\varphi_{(e_i,J)}=\Gamma_{(e_i,J)}\\
    f_{i,(J)}+\varphi_{(J,e_i)}=\Gamma_{(J,e_i)}
    \end{cases}
    \hskip5pt \mbox{for}\ |J|\geq 2,j_i=0;\label{3.2.13}\\
z^I,\ \-{z}^I:&\ \ \  \begin{cases}
    -g_{(I)}+\sum\limits_{i=1}^{n}\left(|z_i|^2f_{i,(I+e_i)}\right)+\varphi_{(I,0)}=\Gamma_{(I,0)}\\
    \sum\limits_{i=1}^{n}\left(|z_i|^2\overline{f_{i,(I+e_i)}}\right)+\varphi_{(0,I)}=\Gamma_{(0,I)}
    \end{cases}
    \hskip5pt \mbox{for}\ |I|\geq 2; \label{3.2.14} \\
z^I\-{z}^J:&\ \ \ \varphi_{(I,J)}=\Gamma_{(I,J)} \ \hbox{for}\ |I|,
|J|\ge 2,\ i_l \cdot j_l=0,\ l=1,\cdots,n. \label{3.2.15}
\end{align}

Here we demonstrate in details  how  the system (\ref{3.2.14}) is uniquely solved. The
others are done similarly  (and, in fact,  more easily). We first
substitute (\ref{2.3}) to (\ref{3.2.14}) and then  collect
coefficients of the zeroth order term, linear terms and higher order terms  in $ v_2,\cdots,v_n$, respectively, while taking  u as a parameter. We
obtain, by Lemma 2.2, the following:

\begin{align}
&\sum_{k}\Gamma_{(I,0)}^{(ke_1)}u^k=-g_{(I)}(u)+2^{1-n}uf_{1,(I+e_1)}+\sum_{i=2}^{n}2^{i-1-n}uf_{i,(I+e_i)}+\sum_{k}\varphi_{(I,0)}^{(ke_1)}u^k;
\label{3.2.17}\\
&\sum_{k}\Gamma_{(I,0)}^{(ke_1+e_j)}u^k=2^{1-j}f_{1,(I+e_1)}+\sum_{i=2}^{j-1}2^{i-1-j}f_{i,(I+e_i)}-2^{-1}f_{j,(I+e_j)}+\sum_{k}\varphi_{(I,0)}^{(ke_1+e_j)}u^k,
\ j\ge 2;\label{3.2.18}\\
 &
\varphi_{(I,0)}^{(k_1e_1+k_2e_2+\cdots+k_ne_n)}=\Gamma_{(I,0)}^{(k_1e_1+k_2e_2+\cdots+k_ne_n)},
k_2+\cdots +k_n\ge 2; \label{3.2.180}
\end{align}
\begin{align}
&\sum_{k}\Gamma_{(0,I)}^{(ke_1)}u^k=2^{1-n}u\overline{f_{1,(I+e_1)}}+\sum_{i=2}^{n}2^{i-1-n}u\overline{f_{i,(I+e_i)}}
+\sum_{k}\varphi_{(0,I)}^{(ke_1)}u^k;\label{3.2.19}\\
&\sum_{k}\Gamma_{(0,I)}^{(ke_1+e_j)}u^k=2^{1-j}\overline{f_{1,(I+e_1)}}+\sum_{i=2}^{j-1}2^{i-1-j}\overline{f_{i,(I+e_i)}}-2^{-1}\overline{f_{j,(I+e_j)}}
+\sum_{k}\varphi_{(0,I)}^{(ke_1+e_j)}u^k,\ j\ge 2;\  \label{3.2.20}\\
&
\varphi_{(0,I)}^{(k_1e_1+k_2e_2+\cdots+k_ne_n)}=\Gamma_{(0,I)}^{(k_1e_1+k_2e_2+\cdots+k_ne_n)},\
k_2+\cdots +k_n\ge 2. \label{3.2.200}
\end{align}

Using the normalization in $\varphi$ and letting $u=0$ in (\ref{3.2.17}) (\ref{3.2.19}), we get
$\Gamma^{(0)}_{(I,0)}=-g_{(I)}(0)+\varphi^{(0)}_{(I,0)}$ and
$\Gamma^{(0)}_{(0,I)}=\varphi^{(0)}_{(0,I)}$. By the normalization
$\varphi^{(0)}_{(I,0)}=\-{\varphi^{(0)}_{(0,I)}}$, we get
$\varphi^{(0)}_{(I,0)}=\-{\Gamma^{(0)}_{(I,0)}}$ and
\begin{align}
g_{(I)}(0)=\-{\Gamma^{(0)}_{(0,I)}}-\Gamma^{(0)}_{(I,0)}.
\end{align}

Sum  up (\ref{3.2.20}) with $j=2,\cdots,n$ and then add it to
(\ref{3.2.19}). By the   the normaliztaion condition
$\varphi_{(0,I)}^{(le_1)}=\varphi_{(I,0)}^{(ke_+e_j)}=0$ for $k\ge 0,
l\ge 1$, we obtain the following:
\begin{align}
f_{1,(I+e_1)}(u)=\sum_{k\ge
1}\-{\Gamma_{(0,I)}^{(ke_1)}}u^{k-1}+\sum_{k\ge
0}\sum_{j=2}^{n}\-{\Gamma_{(0,I)}^{(ke_1+e_j)}}u^k.
\end{align}
Subtracting the complex conjugate of (\ref{3.2.20}) from (\ref{3.2.18}), we obtain
\begin{align}
\varphi_{(I,0)}^{(ke_1+e_j)}=\Gamma_{(I,0)}^{(ke_1+e_j)}-\-{\Gamma_{(0,I)}^{(ke_1+e_j)}},\ j\ge 2,\ k\ge 0.
\end{align}
From (\ref{3.2.17}) and (\ref{3.2.20}), we can similarly get
\begin{align}
g_{(I)}(u)=\sum_{k=0}^{\infty}\left(\-{\Gamma_{(0,I)}^{(ke_1)}}-\Gamma_{(I,0)}^{(ke_1)}\right)
u^k,\ \ |I|\ge 2.
\end{align}

Back to the equation (\ref{3.2.20}), we can inductively get:
\begin{align}
f_{j,(I+e_j)}(u)= \sum_{k\ge
1}\overline{\Gamma_{(0,I)}^{(ke_1)}}u^{k-1}+\sum_{k\ge
0}\left(\sum_{i=0}^{n-j-1}\overline{\Gamma_{(0,I)}^{(ke_1+e_{n-i)}}}-
\overline{\Gamma_{(0,I)}^{(ke_1+e_{j})}}\right)u^k\ \mbox{for}\
2\leq j\leq n.
\end{align}
Similarly, we get from (\ref{3.2.10}) the following
\begin{align}
&g_{(0)}(u)=\sum_{k\ge 2}\left(-\Gamma_{(0,0)}^{(ke_1)}u^k\right)-Re
  \left(\sum_{k\ge 1;\ j=2,\cdots,n}
  \Gamma_{(0,0)}^{(ke_1+e_j)}u^{k+1}\right);\\
& f_{h,(e_h)}(u)=\frac{1}{2}\sum_{k\ge 1}\left(-\sum_{j=2}^{h-1}Re(
  \Gamma_{(0,0)}^{(ke_1+e_j)}u^k)-2Re (
  \Gamma_{(0,0)}^{(ke_1+e_h)}u^k)\right), \ h\ge 2;\\
&\varphi_{(0,0)}=\Gamma_{(0,0)}-\sum\limits_{k\geq
2}\Gamma_{(0,0)}^{(ke_1)}u^k-Re(\sum\limits_{k\geq
1,j=2,\cdots,n}\Gamma_{(0,0)}^{(ke_1+e_j)}u^kv_j).
\end{align}
From (\ref{3.2.11}), we obtain the following:
\begin{align}
& f_{1,(e_1+e_j)}(u)=\sum_{k\ge
1}\-{\Gamma_{(0,e_j)}^{(ke_1)}}u^{k-1}+\sum_{k\ge
0}\sum_{i=2}^{n}\-{\Gamma_{(0,e_j)}^{(ke_1+e_i)}}u^k;
\\
&f_{i,(e_j+e_i)}(u)= \sum_{k\ge
1}\overline{\Gamma_{(0,e_j)}^{(ke_1)}}u^{k-1}+\sum_{k\ge
0}\left(\sum_{l=0}^{n-i-1}\overline{\Gamma_{(0,e_j)}^{(ke_1+e_{n-l)}}}-
\overline{\Gamma_{(0,e_j)}^{(ke_1+e_{i})}}\right)u^k\ \mbox{for}\ 2\leq i\leq n;\\
&g_{(e_j)}(u)=\sum_{k=1}^{\infty}\left(\-{\Gamma_{(0,e_j)}^{(ke_1)}}-\Gamma_{(e_j,0)}^{(ke_1)}\right) u^k;\\
&\varphi_{(e_j,0)}^{(ke_1+e_l)}=\Gamma_{(e_j,0)}^{(ke_1+e_l)}-\-{\Gamma_{(0,e_j)}^{(ke_1+e_l)}},\ l\ge 2,\ k\ge 0;\\
&\varphi_{(e_j,0)}^{(k_1e_1+k_2e_2+\cdots+k_ne_n)}=\Gamma_{(e_j,0)}^{(k_1e_1+k_2e_2+\cdots+k_ne_n)},
\ k_2+\cdots +k_n\ge 2;\\
&\varphi_{(0,e_j)}^{(k_1e_1+k_2e_2+\cdots+k_ne_n)}=\Gamma_{(0,e_j)}^{(k_1e_1+k_2e_2+\cdots+k_ne_n)},
\ k_2+\cdots +k_n\ge 2.
\end{align}
From (\ref{3.2.12}), we get
\begin{align}
&f_{i,(e_j)}(u)=\sum_{k=1}^{\infty}{\Gamma_{(e_j,e_i)}^{(ke_1)}}u^k, \ i<j;\\
&\varphi^{ke_1}_{(e_i,e_j)}=\Gamma^{(ke_1)}_{(e_i,e_j)}-\-{\Gamma^{(ke_1)}_{(e_j,e_i)}}, \ i<j,\ k\geq 1;\\
&\varphi_{(e_i,e_j)}^{(k_1e_1+k_2e_2+\cdots+k_ne_n)}=\Gamma_{(e_i,e_j)}^{(k_1e_1+k_2e_2+\cdots+k_ne_n)}
,\ \rm{for} \ k_2+\cdots+k_n\geq 1.
\end{align}
From (\ref{3.2.13}), we obtain
\begin{align}
&f_{i,(J)}(u)=\sum_{k\ge 0}{\Gamma_{(J,e_i)}^{(ke_1)}}u^k,\ \ 1\le i\le n ,\\
&\varphi^{(ke_1)}_{(e_i,J)}=\Gamma_{(e_i,J)}^{(ke_1)}-\-{\Gamma_{(J,e_i)}^{(ke_1)}},\ \ 1\le i\le n ,k\geq0,\\
&\varphi_{(J,e_i)}^{(k_1e_1+k_2e_2+\cdots+k_ne_n)}=\Gamma_{(J,e_i)}^{(k_1e_1+k_2e_2+\cdots+k_ne_n)},\
\rm{for} \ k_2+\cdots+k_n\geq 1,\\
&\varphi_{(e_i,J)}^{(k_1e_1+k_2e_2+\cdots+k_ne_n)}=\Gamma_{(e_i,J)}^{(k_1e_1+k_2e_2+\cdots+k_ne_n)},\
\rm{for} \ k_2+\cdots+k_n\geq 1.
\end{align}
where $|J|\ge 2$ and $j_i=0$.

Summarizing the solutions just obtained,  we have the following
formula: (One can also directly verify that they are indeed the
solutions of (\ref{2.9})
with the
normalization conditions given in (\ref{2.6}) and (\ref{2.8}))
\begin{equation}
\begin{array}{cll}
F_1(z,u)&=&z_1+f_1(z,u)=z_1+\sum\limits_{k\ge 0,j_1=0,|J|\geq 1
  }z^J\Gamma_{(J,e_1)}^{(ke_1)}u^k+\sum\limits_
  {|I|\geq
  1}z^{I+e_1}S^{(1)}_I,\\
F_h(z,u)&=&z_h+f_h(z,u)=z_h+\frac{1}{2}z_h\sum\limits_{k\ge
1}\big(-\sum\limits_{j=2}^{h-1}Re(
  \Gamma_{(0,0)}^{(ke_1+e_j)}u^k)-2Re (
  \Gamma_{(0,0)}^{(ke_1+e_h)}u^k)\big)\\ &&+\sum\limits_{k\ge 1,i>
  h}z_i\Gamma_{(e_i,e_h)}^{(ke_1)}u^k
+\sum\limits_{k\ge 0,j_h=0,|J|\geq
  2}z^J\Gamma_{(J,e_h)}^{(ke_1)}u^k+\sum\limits_{|I|\geq
  1}z^{I+e_h}S^{(h)}_I,\ \ \  n\ge h\geq 2,\\
G(z,u)&=&u+g(z,u)=u+\left(-\sum\limits_{k\ge
2}\Gamma_{(0,0)}^{(ke_1)}u^k-Re (\sum\limits_{k \geq 1,j=2,\cdots,n}
  \Gamma_{(0,0)}^{(ke_1+e_j)}u^{k+1})\right)\\ &&+\sum\limits_{k\ge 0,|I|\geq
  1}z^Iu^k\left(\overline{\Gamma_{(0,I)}^{(ke_1)}}-\Gamma_{(I,0)}^{(ke_1)}
 \right),\\
\varphi&=&\Gamma(z,\bar{z})+g(z,u)-2Re\left(\sum\limits_{i=1}^{n}\left(\overline{z_i}f_i(z,u)\right)\right),
\label{2.10}
\end{array}
\end{equation}
where
\begin{equation}\left\{
\begin{array}{l}
S^{(1)}_I=\sum\limits_{k\ge 1}\overline{\Gamma_{(0,I)}^{(ke_1)}}u^{k-1}+\sum\limits_{k\ge 0}\sum\limits_{i=2}^{n}\overline{\Gamma_{(0,I)}^{(ke_1+e_{i})}}u^{k},\\
S^{(h)}_I=\sum\limits_{k\ge
1}\overline{\Gamma_{(0,I)}^{(ke_1)}}u^{k-1}+\left(\sum\limits_{k\ge
0}\sum\limits_{i=0}^{n-h-1}\overline{\Gamma_{(0,I)}^{(ke_1+e_{n-i)}}}u^{k}\right)-\sum\limits_{k\ge
0} \overline{\Gamma_{(0,I)}^{(ke_1+e_{h})}}u^k,\ \mbox{for}\ 2\leq
h\leq n.
\end{array}
\right. \label{2.11}
\end{equation}
This completes the proof of Theorem 2.3. $\endpf$
\\

Let $(M,0)$ be as in (\ref{2.1}). We say that $(M^*,0)$ is a formal
pseudo-normal form for $(M,0)$ if $(M^*,0)$ is formally equivalent
to $(M,0)$ and $M^*$ is defined by $w=|z|^2+\varphi$ with $\varphi$
satisfying the normalizations in (\ref{2.9}). We notice that
pseudo-normal forms of $(M,0)$ are not unique. Furthermore, we have
the following observations:

\medskip

{\bf Remark 2.4}: {\bf (A)}.The pseudo-normal form obtained in
Theorem 2.3 contains information reflecting both the singular CR
structure and partial strongly pseudoconvex CR structure at the
point under study. For instance, the following submanifold in
${\CC}^3$ is given in a pseudo-normal form:
\begin{equation}
M:\ w=|z|^2+2Re\sum_{j_1+ j_2\ge
3}\left( a_{j_1j_2}z_1^{j_1}z_2^{j_2}\right)+\sum_{j_1\ge 2, j_2\ge
2}b_{j_1 \-{j_2}}z_1^{j_1}\-{z_2}^{j_2}.
\label{2.hu-01}
\end{equation}
 Here the  harmonic terms
$Re\sum_{j_1+ j_2\ge 3}\left( a_{j_1j_2}z_1^{j_1}z_2^{j_2}\right)$
are presented  due to the nature of CR singularity of $M$ at $0$,
which may be compared with the Moser pseudo-normal form in [Mos] in
the pure CR singularity setting. Typical mixed terms like $\sum_{j_1\ge 2,
j_2\ge 2}b_{j_1 \-{j_2}}z_1^{j_1}\-{z_2}^{j_2}$ are associated with  the partial  CR
structure near $0$, which can be compared with the Chern-Moser
normal form in the pure CR setting [CM].

{\bf (B)}. Suppose that $M$
 is defined by a formal equation of the form:
$w=|z|^2+E(z,\-{z})$ with $Ord (E)\ge 3$ and
$\-{E(z,\-{z})}=E(z,\-{z})$. In the normalized map
$H(z,w)=(F(z,w),G(z,w))$ transforming $M$ into its normal form in
Theorem 2.3,  the $w$-component $G(z,u)$ is only a function in $u$
and is formally real-valued, by the formula in (\ref{2.10}). This is
due to the fact that  the $\Gamma$ in (\ref{2.10}) obtained from
each induction stage in the process of the proof of Theorem 2.3 is
formally real-valued. Hence, the $\varphi$ in the
pseudo-normalization of $M$ obtained in Theorem 2.3 is also formally
real-valued. However, fundamentally different from the two
dimensional case, this is no longer true for a general $M$. Indeed,
we will see in Theorem 3.5 that $M$ can be formally flattened if and
only if its pseudo-normal form is given by a formal real-valued
function.

\bigskip

\section{ Normalization of  holomorphic maps by   automorphisms  of the  quadric }

In this section, we first compute the isotropic automorphism group
of the model space $M_{\infty}\subset \mathbb{C}^{n+1}$ defined by
the equation:
$w=\sum_{i=1}^n|z_i|^2$.
Write $Aut_0(M_\infty)$ for the set of biholomorphic self-maps of
$(M_\infty, 0)$. We have the following:

\medskip
{\bf Proposition 3.1:}\hspace{0.2cm} {\it $Aut_0(M_\infty)$ consists of
 the transformations given in the following (\ref{3.1}) or (\ref{3.2}) :\\
\begin{equation}
 \left\{
\begin{array}{l}
z'=b(w)\frac
          {wa(w)-\frac{\langle z,\bar{a}(w)\rangle}{\langle
                a(w),\bar{a}(w)\rangle}a(w)
                +\sqrt{1-wa(w)\bar{a}(w)}\left(z-\frac{\langle z,\bar{a}(w)\rangle}{\langle a(w),\bar{a}(w)\rangle}a(w)\right)}
          {1-\langle z,\bar{a}(w)\rangle}
  U(w)\\
w'=b(w)\bar{b}(w)w
\end{array}
\right. \label{3.1}
\end{equation}
\begin{equation}
(z',w')=\left(b(w)zU(w),b(w)\bar{b}(w)w\right). \label{3.2}
\end{equation}
where $a=(a_1,\cdots,a_n)$,
$\sum_{j=1}^na_j(0)\bar{a_j}(0)<1$,
 $\langle z,\bar{a} \rangle
=\sum\limits_{i=1}^n\bar{a_i}z_i$, $b(0)\neq 0$, $a(0)\neq 0$,
$U(Re(w))$ is a unitary matrix and $a(w),b(w),U(w)$ are holomorphic
in $w$.}

\medskip
{\it Proof of Proposition 3.1:} Write $w=x+\sqrt{-1}y$. Let
$(F,G)\in Aut_0(M_\infty)$. Then $Im(G(z,|z|^2))\equiv 0$ for
$z\approx 0$. Since $M_\infty$ bounds a family of balls near 0
defined by
$$
B_r=\left\{(z,w)\in \mathbb{C}^{n+1}:w=x+\sqrt{-1}y,y=0,x=r^2 \geq
|z^2| \right\}.
$$
We see that $Im(G(z,x))\equiv 0$ for $z\approx 0$ and $x(\in
{\RR})\approx 0$. Therefore, $G(z,w)=G(w)=cw+o(w)$ $(c>0)$ is
independent of $z$ and takes  real value when $w=x$ is real.
Now  $F(z,r^2)$ must be a biholomorphic map from $|z|^2 < r^2$ to
$|z|^2 <G(r^2)$ for any $r>0$. Using the explicit expression for
automorphisms of the unit ball (see [Rud]), we obtain either:
\begin{equation}
F(z,r^2)=\sqrt{G(r^2)}\frac
          {a(r)-\frac{\langle {z \over r},\bar{a}(r)\rangle}{\langle
                a(r),\bar{a}(r)\rangle}a(r)
                +v\left({z \over r}-\frac{\langle {z \over r},\bar{a}(r)\rangle}{\langle a(r),\bar{a}(r)\rangle}a(r)\right)}
          {1-\langle {z \over r},\bar{a}(r)\rangle}
  U(r)\\
  \label{3.3}
\end{equation}
where $U(r)$ is a unitary  matrix and $v=\sqrt{1-a(r)\bar{a}(r)},\
a\not =0$, or we have
\begin{equation}
F(z,r^2)=\sqrt{G(r^2)}({z \over r})U(r).\label{3.3'}
\end{equation}
Write $G(x)=xb(x)\-{b}(x)$ with $b(0)\not = 0$  and $b(w)$
holomorphic in $w$. In the case of (\ref{3.3'}),
$F(z,x)=b(x)zU(r)e^{\sqrt{-1}\theta(x)}$ is real analytic, where
$\theta(x)$ is real-valued real analytic function in $x$ . Hence,
$b(x)U(r)e^{\sqrt{-1}\theta(x)}$ is the Jacobian matrix of $F$ in
$z$. Since both $e^{\sqrt{-1}\theta(x)}$ and $b(x)(\not =0)$ are
real analytic for $x\approx 0$, we conclude that $U(r)$ is real
analytic in $x$. Hence, $U(w)$ is also holomorphic in $w$. Still
write $U(x)$ for $U(x)e^{i\theta(x)}$. We see the proof of
Proposition 3.1 in the case of (\ref{3.2}).

 Suppose that $a \neq 0$. Still write $G(w)=wb(w)\-{b}(w)$ with $b(0)\not = 0$. We  have
$$
F(z,r^2)=b(r^2)\frac
          {ra(r)-\frac{\langle z,\bar{a}(r)\rangle}{\langle
                a(r),\bar{a}(r)\rangle}a(r)
                +v\left(z-\frac{\langle z,\bar{a}(w)\rangle}{\langle a(r),\bar{a}(r)\rangle}a(r)\right)}
          {1-\langle z,{\bar{a}(r) \over r} \rangle}e^{i \theta}
  U(r)\\
$$
Since $f(z,w)$ is holomorphic in $(z,w)$ and
$f(0,w)=b(w)\sqrt{w}a(\sqrt{w})U^*(\sqrt{w})$ with
$U^*=e^{i\theta}U$, we see that $\sqrt{w}a(\sqrt{w})U^*(\sqrt{w})$
is holomorphic in $w$ . In particular, $|a(\sqrt{w})|^2$ is real
analytic in $w$. Moreover,
$$
\frac{\partial F}{\partial
z_i}(0,w)=b(w)\left(\frac{|a|^2-v-1}{|a|^2}\bar{a_i}a+ve_i\right)U^*(\sqrt{w})
$$
is analytic. Since $\sqrt{w}a(\sqrt{w})U^*(\sqrt{w})$ is real
analytic, we see that
$$
\left(\frac{|a|^2-v-1}{|a|^2}\bar{a_i}a+ve_i\right)U^*(\sqrt{w})\overline{U^*(\sqrt{w})}^t\overline{a(\sqrt{w})}^tr
=\left((|a|^2-v-1)+v\right)r\bar{a_i}
$$
is real analytic, too.  Here $(\cdot)^t$ denotes the matrix
transpose. Since $(|a|^2-v-1)+v=|a|^2-1$ is real analytic, we
conclude that both $ra_i$ and $a_i/r$ are real analytic in $w$. Since both
$\sqrt{w}a(\sqrt{w})U^*(\sqrt{w})$ and $ra_i$ are real
analytic, we  see that $U^*(\sqrt{w})$ is  real analytic.
Still denote $a$ for $a/r$,  we further obtain the following with the given properties stated in the Proposition:
\begin{equation*}
 \left\{
\begin{array}{l}
F(z,w)=b(w)\frac
          {wa(w)-\frac{\langle z,\bar{a}(w)\rangle}{\langle
                a(w),\bar{a}(w)\rangle}a(w)
                +\sqrt{1-wa(w)\bar{a}(w)}\left(z-\frac{\langle z,\bar{a}(w)\rangle}{\langle a(w),\bar{a}(w)\rangle}a(w)\right)}
          {1-\langle z,\bar{a}(w)\rangle}
  U^*(w)\\
G(w)=b(w)\bar{b}(w)w.
\end{array}
\right.
\end{equation*}
This completes the proof of Proposition 3.1. $\endpf$

\medskip

{\bf Remark 3.2}: In  Proposition 3.1, if we let $a(w), b(w), U(w)$
be formal power series in $w$ with $a(0),\ b(0)\not =0$ and $\langle
a(0),\-{a}(0)\rangle<1$, $U(x)\cdot U(x)^t=I$, then (\ref{3.1}) and
(\ref{3.2})  give formal automorphisms of $M_\infty$, which are not
convergent. Write the set of automorphisms obtained in this way as
$aut_0(M_\infty)$. One may prove that $aut_0(M_\infty)$ consists of
all the formal automorphisms of $(M_\infty, 0)$ .

\medskip

 We now suppose
that $H=(F,G)$ is a formal equivalence self-map of
$({\CC}^{n+1},0)$, mapping a formal submanifold  of the form
$w=|z|^2+O(|z|^3)$ to  a submanifold of the form $w=|z|^2+O(|z|^3)$.
The following lemma shows that we can always normalize $H$ by
composing it from the left with an element from $aut_0(M_\infty)$ to
get a normalized mapping. This fact will be used in the proof of
Theorem 1.
 In what
follows, we
set $v(g,a)=\sqrt{1-g \cdot a(g)\cdot\bar{a}(g)}$.
\bigskip

{\bf Lemma 3.3:}\hspace{0.2cm} {\it There exists a unique
automorphism $T\in aut_0(M_\infty)$   such that $T\circ H$ satisfies
the normalized condition in (\ref{2.6}). When $H$ is biholomorphic,
$T\in Aut_0(M_\infty)$ }
\bigskip

{\it Proof of Lemma 3.3:} First, it is easy to see that by composing
an automorphism of the form $w'=|c|^2w,\ z'=czU$, we can assume that
$F=z+O_{wt}(2)$ and $G=w+O_{wt}(3)$ (see [Hu1]). Here $c$ is a non-zero constant
and $U$ is a certain $n\times n$-unitary matrix.

 Let $
b(w)=1,a_j=\alpha_j(w),a_1=\cdots=a_{j-1}=a_{j+1}=\cdots=a_n=0$, and
$U=I$ in (\ref{3.1}). We get the following automorphism of
$M_\infty$
$$
T_j=\left(\frac{v(w,\alpha_j)z_1}{1-\bar{\alpha_j}z_j},\cdots,\frac{v(w,\alpha_j)z_{j-1}}{1-\bar{\alpha_j}z_j},\frac{z_j-w\alpha_j}{1-\bar{\alpha_j}z_j}
   ,\frac{v(w,\alpha_j)z_{j+1}}{1-\bar{\alpha_j}z_j},\cdots,\frac{v(w,\alpha_j)z_n}{1-\bar{\alpha_j}z_j},w\right).
$$
Write
\begin{equation}
\left\{
\begin{array}{l}
H_j=(_{(j)}F,_{(j)}G)=T_j  \circ T_{j-1} \circ  \cdots \circ T_1
\circ H,\
H_0=H;\\
\alpha_j=\frac{_{(j-1)}F_{j,(0)}(u)}{_{(j-1)}G_{(0)}(u)}\circ
\left(_{(j-1)}G_{(0)}(u)\right)^{-1}.
\end{array}
\right.
\end{equation}
Then a direct computation shows that $(_{(j)}F)_{i,(0)}(u)=0$ for
$1\leq i \leq j$. In particular, we have
 $(_{(j)}F)_{i,(0)}(u)=0$ for all $1\leq i\leq n$.\\
Still write $H$ for $H_n$. Next, for  $i<j$, let $ b(w)=1, a=0$, and
let
\begin{equation*}
U_j^i= \left (
\begin{array}{ccccc}
I \ \ &0\  \ &0\ \ &0\ \ &0\\
0 &\cos(\theta_j^i) &0 &-\sin(\theta_j^i) &0\\
0 &0 &I &0 &0\\
0 &\sin(\theta_j^i) &0 &\cos(\theta_j^i) &0\\
0 &0 &0 &0 &I\\
\end{array}
\right )
\end{equation*}
in (\ref{3.1}), where $\cos(\theta_j^i)$ is at the $i^{th}$ row  and
the $j^{th}$ column. Then we get an automorphism $T_{j}^{i}$. Set
\begin{equation}
\begin{array}{l}
H_{j}^{i}=(_{j}^{i}F,_{j}^{i}G)=T_j^i \circ \cdots \circ T_{i+1}^i
\circ T_n^{i-1} \circ
\cdots T_i^{i-1} \circ \cdots \circ T_n^1 \circ \cdots \circ T_2^1 \circ H,\\
\theta_j^i=\left\{
\begin{array}{ll}
\tan^{-1}\left(\frac{(_n^{i-1}F)_{j,(e_i)}}{(_n^{i-1}F)_{i,(e_i)}}\right)
 \circ \left(_n^{i-1}G_{(0)}(w)\right)^{-1},
 \hskip 0.5cm &j=i+1;\\
\tan^{-1}\left(\frac{(_{j-1}^{i}F)_{j,(e_i)}}{(_{j-1}^{i}F)_{i,(e_i)}}\right)
\circ \left(_{j-1}^{i}G_{(0)}(w)\right)^{-1},
&j\neq i+1.\\
\end{array}
\right.
\end{array}
\end{equation}
Then we can inductively prove that $H_j^i$ satisfies
$$
(_j^iF)_{(0)}=0\ , \ (_j^iF)_{k,(e_l)}=0\ \mbox{for}\ l=i,i+1 \leq
k\leq j\ \mbox{or} \ l<i,l+1\leq k \leq n.
$$
In particular, we see that $H_n^{n-1}$ satisfies
$(_n^{n-1}F)_{(0)}=0\ ,
\ (_n^{n-1}F)_{i,(e_j)}=0\ \mbox{for}\ 1 \leq j < i \leq n$.\\
Still write $H$ for $H_n^{n-1}$ and set $H'=T \circ H=(F',G')$ with
$$
T=\left(d(w)z,d(w)\bar{d}(w)w\right)\ ,\
d=\frac{1}{F_{1,(e_1)}\left(w\right)}
\circ\left(G_{(0)}(w)\right)^{-1} .
$$
Then $H'$ satisfies
$$
(F')_{(0)}=0\ ,\ (F')_{1,(e_1)}=1\ , \ (F')_{i,(e_j)}=0\ \mbox{for}\
1 \leq j < i \leq n.
$$
At last, a composition from the left with the rotation map as
follows:
\begin{equation*}
\hat{T}=(z_1,\beta_2z_2,\cdots,\beta_nz_n,w),\
\beta_i=\frac{(\bar{F'})_{i,(e_i)}(w)}{\sqrt{(F')_{i,(e_i)}(w)\cdot
(\bar{F'})_{i,(e_i)}(w)}}
 \circ \left({G'}_{(0)}(w)\right)^{-1}
\end{equation*}
makes $H'$ satisfy the normalization condition (\ref{2.6}). This
proves the existence part of the lemma.

Next, suppose both $H=(F,G)=(z+O_{wt}(2),w+O_{wt}(3))$
 and
$\hat{H}=(\hat{F},\hat{G})= T\circ H=(z+O_{wt}(2),w+O_{wt}(3))$
satisfy the normalization condition (\ref{2.6}). Here T is an
automorphism of $M_\infty$.
Then $T$ must be of the form in (\ref{3.2}), for $T(0,w)=0$. Hence,
$$
T=(b(w)zU(w),b(w)\bar{b}(w)w).
$$
By the normalization condition (\ref{2.6}) on $H,\ \hat{H}$, we have
\begin{equation}
 \left(
\begin{array}{cc}
  \begin{array}{ll}
       1&\\&\hat{F}_{2,(e_2)}
  \end{array}&0\\ \ast&
  \begin{array}{ll}
       \ddots&\\&\hat{F}_{n,(e_n)}
  \end{array}
\end{array} \right)
=b\left(G_{(0)}(w)\right)
 \left(
\begin{array}{cc}
    \begin{array}{ll}
       1&\\&F_{2,(e_2)}
    \end{array}&0\\ \ast&
    \begin{array}{ll}
       \ddots&\\&F_{n,(e_n)}
    \end{array}
\end{array} \right)U(G_0(w)).
\end{equation}
with $U(x)$ unitary and
$Im(\hat{F}_{i,(e_i)}(0,u))=Im(F_{i,(e_i)}(0,u))=0$. Considering the
norm of the first row of the right hand   side, we get
$b(G_{(0)}(w))\cdot \-{b}(G_{(0)}(w))=1$ in case
$G_{(0)}(w)=\-{G_{(0)}(w)}$. Since $G_0(w)=w+o(w)$, this implies
that $b(w)\-{b}(w)\equiv 1$ and thus $T=(b(w)zU(w),w)$. Write
$$
b(w) U(w)=\widetilde{U}(w)=\left(
  \begin{array}{ccc}
   u_{11}&\cdots&u_{nn}\\
   \vdots&\ddots&\vdots\\
    u_{n1}&\cdots&u_{nn}.
   \end{array}\right).
$$
We  notice that  $\widetilde{U}$ is a lower triangular matrix and is
unitary when $w=x$. Thus
we have $u_{ii}(w)\-{u_{ii}}(w)= 1$  and $u_{ij}=0$ for $i\not = j$.
Notice that
$$
u_{11}\equiv 1,\ \hat F_{i,(e_i)}(w)=u_{ii}(w)\cdot F_{i,(e_i)}(w)\
\rm{for} \ 2\leq i\leq n.
$$
Since  $\hat F_{i,(e_i)}(x),F_{i,(e_i)}(x)=1+o(x)$ are real, we get
$u_{ii}(x)=1$.
This proves the uniqueness part of the lemma. $\endpf$

\bigskip

{\bf Lemma 3.4:}\hspace{0.2cm} {\it Suppose that $H$ with $H(0)=0$
is an equivalence map from $w=|z|^2+\varphi(z,\bar{z})$ to
$w'=|z'|^2+\varphi'(z',\bar{z'})$. Here $\varphi$ and $\varphi'$ are
normalized as in (\ref{2.8}). Let $s,s'$ be the lowest order of
vanishing in $\varphi$ and $\varphi'$, respectively, then $s=s'$.
}

{\it Proof of Lemma 3.4}: We seek for a contradiction if $s\not =
s'$.  Assume, for instance,  that $s<s'$.
 Let $T$ be an automorphism of $M_\infty$  with $T \circ H$ being
normalized as in (\ref{2.6}). Suppose that $T$ transforms
$w'=|z'|^2+\varphi'(z',\bar{z'})$ to
$w''=|z''|^2+\varphi''(z'',\overline{z''})$ with $s''$ the lowest
vanishing order for $\varphi''$. We claim that $s'=s''$. Suppose
not. We assume, without loss of generality,  that $s'<s''$.
Write the linear part  of $T$ (in $(z',w')$) as
$(z''=z'B+Dw',w''=dw')$ with $B\in GL(n,{\CC}),\ d\not =0$. Then a
direct computation shows that
$$\varphi''^{(s')}(z'B,\-{z'B})=d\cdot \varphi'^{(s')}(z',z').$$
This is a contradiction.


Now, $T \circ H$ transforms $w=|z|^2+\varphi$ to
$w''=|z''|^2+\varphi''$ with $T\circ H,\ \varphi$ being normalized
as in (\ref{2.6}) and (\ref{2.8}), respectively. Also  $s<s''$. we
see that $T\circ H$ transforms $w=|z|^2+\varphi^{(s)}$ to $w=|z|^2$,
moduling $O(|(z_1,\cdots,z_n)|^{s+1})$. This contradicts the
uniqueness part of Theorem 2.3. The proof of Lemma 3.4 is complete.
$\endpf$

\bigskip
We say that a formal submanifold $(M,0)$ of real dimension $2n$
defined by (\ref{2.1}) can be formally flattened if there is a
formal change of coordinates $(z',w')=H(z,w)$ with $H(0)=0$ such
that in the new coordinates $(M,0)$ is defined by a formal function
of the form $w'=E^*(z',\-{z'})$ with
$E^*(z',\-{z'})=\-{{E^*(z',\-{z'})}}$. We also say a pseudo-normal
form of $(M,0)$ given by $w=|z|^2+\varphi(z,\-{z})$ with $\varphi$
satisfying the normalizations in (\ref{2.8}) is a flat pseudo-normal
form if $\varphi$ is formally real-valued. An immediate application
of Lemma 3.3 and Remark 2.4 (b) is that if $(M,0)$ has a flat
pseudo-normal form, then all of its other pseudo-normal forms are
flat. Indeed, for a given pseudo-normal form of $(M,0)$, there is a
formal equivalence map $H$ mapping it into $Imw=0$. Now, by Lemma
3.3, we can compose $H$ with an element $T$ of $aut_0(M_\infty)$ to
normalize $H$. Next, since $T$ maps any flattened submanifold to a
flattened submanifold, there is a formal transformation $H^*$ such
that $H^*\circ T\circ H$ maps the pseudo-normal form given at the
beginning to a flat one. On the other hand, since $H^*\circ T\circ
H$ satisfies the normalizations in (\ref{2.6}), by Theorem 2.3, we
see that $H^*\circ T\circ H=id$ and two pseudo-normal forms are the
same. Summarizing the above, we proved the following:

\medskip
{\bf Theorem 3.5}: {\it Let $(M,0)$ be a formal submanifold defined by an equation of the form: $w=|z|^2+E(z,\-{z})$ with $E=O(|z|^3)$.
Then the following statements are equivalent:

\noindent
(I).  $(M,0)$ can
be flattened

\noindent (II). $(M,0)$ has a flat pseudo-normal form. Namely, $M$
has a pseudo-normal form given by an equation of the form:
$w'=|z'|^2+\varphi(z',\-{z'})$ with $\varphi$ satisfying the
normalizations in (\ref{2.8}) and the reality condition
$\varphi(z',\-{z'})=\-{\varphi(z',\-{z'})}.$

\noindent (III). Any pseudo-normal form of $(M,0)$ is flat.}

\bigskip
{\bf Remark 3.6}: By Theorem 3.5, we see that $M$ defined in
(\ref{2.hu-01}) can be formally flattened if and only if
$b_{i\-{j}}=\-{b_{j\-{i}}}$ for all $i,j.$

\bigskip

\section{Proof of Theorem 1}

We now give a proof of Theorem 1 by using the rapidly convergent power series method. We let
 $M$ be defined by \\
\begin{equation}
w=\Phi(z,\bar{z})=|z|^2+E(z,\bar{z})
\label{4.1}
\end{equation}
where $E(z,\xi)$ is holomorphic near $z=\xi=0$ with vanishing order
$\geq 3$. Assume that $H=(F,G)=(z+f,w+g)$ is a   formal map
 satisfying the normalization condition
in (\ref{2.2}). We define
\begin{equation}
R=(r_1,r_2,\cdots,r_n)=(2^{-\frac{n-2}{2}}r,2^{-\frac{n-2}{2}}r,2^{-\frac{n-3}{2}}r\cdots,2^{-\frac{1}{2}}r,r).
\label{4.2}
\end{equation}
Then
$|R|^2=2^{-(n-2)}r^2+\sum_{h=2}^{n}(2^{-\frac{n-i}{2}}r)^2=2r^2.$
Define the domains:\\
\begin{equation}
\begin{array}{l}
\Delta_r=\{(z,w):|z_i|<r_i,|w|<2r^2\},\\
D_r=\{(z,\xi):|z_i|<r_i,|\xi_i|<r_i\ \mbox{for}\ 1\leq i \leq n\}.\\
\end{array}
\label{4.3}
\end{equation}
When $E(z,\xi)$ is defined over $\-{D_r}$, we set the norm of
$E(z,\bar{z})$ on $D_r$ by
\begin{equation}
\|E\|_r=\sup\limits_{(z,\xi) \in D_r}|E(z,\xi)|. \label{4.4}
\end{equation}
Also for a holomorphic map $h(z,w)$ defined on $\-{\Delta_r}$, we
define
\begin{equation}
|h|_r=\sup\limits_{(z,w) \in \Delta_r}|h(z,w)|. \label{2.4.5}
\end{equation}
After a scaling transformation $(z,\xi,w)\longrightarrow
(az,a\xi,a^2w)$, we may assume that $E$ is holomorphic on $\-{D_1}$
with $|E|_1 \leq \eta$ for a given small $\eta>0$.

Suppose that $H$ maps $M$ to the quadric $w'=|z'|^2$. Then we have
the following equation:
\begin{equation}
\begin{array}{l}
E(z,\bar{z})+g(z,\Phi)
=2Re\Big(\sum\limits_{i=1}^{n}\bar{z_i}f_i(z,\Phi)\Big)
+|f(z,\Phi)|^2.
 \label{4.5}
\end{array}
\end{equation}
We consider the following  linearized equation of (\ref{4.5}) with $(f,g,\varphi)$ as its unknowns:
\begin{equation}
\begin{array}{l}
E(z,\bar{z})=-g(z,u)+2Re\Big(\sum\limits_{i=1}^{n}\bar{z_i}f_i(z,u)\Big)+\varphi(z,\bar{z}),
\end{array}
\label{4.6}
\end{equation}
where $\varphi$ satisfies (\ref{2.8}). The unique solution of
(\ref{4.6}) is  given in the formula (\ref{2.10}).
However,  we  will make  a certain truncation to $(f,g,\varphi)$ to
faciliate the estimates.
 Suppose that $Ord(E)\geq d\ge 3$. Set
\begin{equation}
\left\{\begin{array}{l}
f=\hat{f}+ O_{wt}(2d-3),\ deg_{wt}(\hat{f})\leq 2d-4,\\
g=\hat{g}+O_{wt}(2d-2),\ deg_{wt}(\hat{g})\leq 2d-3.
\end{array}\right.
\label{4.7}
\end{equation}
Define
\begin{equation*}
\hat{F}=z+\hat{f}\ , \ \hat{G}=w+\hat{g},\
\hat{H}=(\hat{F},\hat{G}).
\end{equation*}
Write $\Theta=(\hat{F},\hat{G})$ and write
\begin{equation}
\begin{array}{l}
\hat{\varphi}(z,\bar{z})=E(z,\bar{z})+\hat{g}(z,u)-2Re\Big(\sum\limits_{i=1}^{n}\bar{z_i}\hat{f}_i(z,u)\Big).
\end{array}
\label{4.6.0}
\end{equation}
Then $\hat{\varphi}(z,\bar{z})-\varphi(z,\bar{z})=O(|z|^{2d-2}).$
 Assume that $M'=\Theta(M)$ is defined by
$w'=|z'|^2+E'(z',\bar{z'})$. Choose $r',\sigma,\varrho,r$ to be such
that
$$
\frac{1}{2}<r'<\sigma<\varrho<r \leq 1,\ \varrho=\frac{1}{3}(2r'+r),
\ \sigma=\frac{1}{3}(2r'+\varrho).
$$
As in the paper of Moser [Mos],  the following lemma will be
fundamental for applying the rapid iteration procedure of Moser to
prove Theorem 1.

\medskip

{\bf Lemma 4.1:}\hspace{0.2cm} {\it Let $M: w=|z|^2+E(z,\bar{z})$ be
as in Theorem 1. Suppose that $Ord(E) \geq d$. Let $\hat{H}$ and
$E'$ be defined above. Then
$
Ord(E') \geq 2d-2.
$}

\medskip
{\it Proof of Lemma 4.1:}  Making use of (\ref{4.6.0}), we have
\begin{equation}
\begin{array}{l}
E'(z',\-{z'})=\Big(\hat{g}(z,\Phi)-\hat{g}(z,u)\Big)-2Re\Big(\sum\limits_{i=1}^{n}\overline{z_i}(\hat{f}_i(z,\Phi)
   -\hat{f}_i(z,u))\Big)-|\hat{f}(z,\Phi)|^2+\hat\varphi(z,\bar{z}).
\end{array}
\label{4.8}
\end{equation}
Since  $Ord(E)\geq d$, by (\ref{2.10}) and (\ref{2.11}),  we see
that $Ord(\hat{f})\geq d-1$ and $Ord(\hat{g}) \geq d$. Hence, we
have
\begin{equation*}
\begin{array}{l}
Ord\Big(\hat{g}(z,\Phi)-\hat{g}(z,u)\Big) \geq  \min\{(d-1)+d,2d-2\}=2d-2,\\
Ord\Big(\hat{f}_i(z,\Phi)-\hat{f}_i(z,u)\Big) \geq \min\{(d-2)+d,2d-3\}=2d-3,\\
Ord\Big(\left|\hat{f}(z,\Phi)\right|^2\Big) \geq 2(d-1)=2d-2.
\end{array}
\end{equation*}
Thus $Ord(E'-\hat{\varphi}) \geq 2d-2$. By  Lemma 3.3 and the
assumption that $w=|z|^2+E$ is formally equivalent to $w=|z|^2$, we
have $s=\infty$. Hence we have $Ord(\varphi) \geq 2d-2$. The lemma
follows.$\endpf$

\medskip
Before proceeding to the estimates of the solution given in (\ref{2.10}), we need
the following lemma:
\bigskip

{\bf Lemma 4.2:} \hspace{0.2cm}{\it If $E$ is holomorphic in
$\-{D_r}$, then we have
\begin{equation*}
\begin{array}{l}
|E^{(ke_1)}_{(I,T)}|\leq
\frac{(k+2)^n\|E\|_r}{R^{I+T}\cdot(2r^2)^k}\ ,\
|E^{(ke_1+e_j)}_{(I,T)}| \leq
\frac{2^n(k+2)^n\|E\|_r}{R^{I+T}(2r^2)^{k+1}}.
 \end{array}
\end{equation*}
}

{\it Proof of Lemma 4.2:} We here give the estimates for
$|E^{(ke_1)}_{(0,I)}|$, $|E^{(ke_1+e_j)}_{(0,I)}|$. The others can
be done similarly. Suppose that $E=\sum a_{i_1\cdots i_nj_1\cdots
j_n}z_1^{i_1} \cdots z_n^{i_n} \overline{z_1}^{j_l}\cdots
\overline{z_1}^{j_l}$. Then by (\ref{2.4}), we have
\begin{equation*}
\begin{array}{lll}
E_{(0,I)}&=&\sum\limits_J a_{j_1\cdots j_n(i_1+j_1)\cdots
(i_n+j_n)}|z_1|^{2j_1}\cdots |z_n|^{2j_n}\\
    &=&\sum\limits_J
    a_{J(I+J)}\left(2^{1-n}(u+\sum\limits_{i=2}^{n}2^{n-i}v_i)\right)^{j_1}
    \Pi_{h=2}^{n}\left(2^{h-n-1}(u+\sum\limits_{i=h+1}^{n}2^{n-i}v_i-v_h)\right)^{j_h}\\
    &=&\sum\limits_J
    a_{J(I+J)}2^{-\left((n-1)j_1+\sum\limits_{h=2}^{n}(n-h+1)j_h\right)}
    \Big(u^{|J|}+\sum\limits_{k=2}^{n}2^{n-k}\left(\sum\limits_{h=1}^{k-1}j_h-j_k\right)u^{|J|-1}v_k\\
    &&\hskip 10pt +O\left(\left|(v_2,\cdots,v_n)\right|^2\right)\Big).\\
\end{array}
\end{equation*}
Thus we obtain
\begin{equation}
\begin{array}{l}
E_{(0,I)}^{(ke_1)}=\sum\limits_{|J|=k}a_{J(I+J)}2^{-\left((n-1)j_1+\sum\limits_{h=2}^{n}(n-h+1)j_h\right)},\\
E_{(0,I)}^{(ke_1+e_l)}=\sum\limits_{|J|=k+1}a_{J(I+J)}2^{-\left((n-1)j_1+\sum\limits_{h=2}^{n}(n-h+1)j_h\right)}2^{n-l}\left(\sum\limits_{h=1}^{l-1}j_h-j_l\right).
\end{array}
\label{4.9}
\end{equation}
By the Cauchy estimates, we get
\begin{equation*}
\begin{array}{lll}
|E_{(0,I)}^{ke_1}|&=&|\sum\limits_{|J|=k}a_{J(I+J)}2^{-\left((n-1)j_1+\sum\limits_{h=2}^{n}(n-h+1)j_h\right)}|\\
     &\leq& \sum\limits_{|J|=k}\frac{\|E\|_r}{R^{I+2J}}2^{-\left((n-1)j_1+\sum\limits_{h=2}^{n}(n-h+1)j_h\right)}\\
     &=& \sum\limits_{|J|=k}\frac{\|E\|_r}{R^{I}}\frac{2^{-\left((n-1)j_1+(n-1)j_2+\cdots+j_n\right)}}
          {(2^{2-n}r^2)^{j_1}\cdot (2^{2-n}r^2)^{j_2} \cdots
          (r^2)^{j_n}}\\
     &\leq& \frac{(k+1)^n\|E\|_r}{R^{I}\cdot(2r^2)^k},\\
|E_{(0,I)}^{ke_1+e_l}|&=&|\sum\limits_{|J|=k+1}a_{J(I+J)}2^{-\left((n-1)j_1+\sum\limits_{h=2}^{n}(n-h+1)j_h\right)}2^{n-l}
\left(\sum\limits_{h=1}^{l-1}j_h-j_l\right)|\\
     &\leq& \sum\limits_{|J|=k+1}\frac{\|E\|_r}{R^{I+2J}}2^{-\left((n-1)j_1+\sum\limits_{h=2}^{n}(n-h+1)j_h\right)}2^{n-l}|\sum\limits_{h=1}^{l-1}j_h-j_l|\\
     &\leq& \sum\limits_{|J|=k+1}\frac{\|E\|_r}{R^{I}\cdot
        (2r^2)^{k+1}}2^n(k+1)\\
     &=&\frac{2^n(k+2)^n\|E\|_r}{R^I\cdot (2r^2)^{k+1}}.
\end{array}
\end{equation*}
Here we have used the fact that
$$
\sharp\left\{(j_1,j_2,\cdots,j_n)\in \mathbb{Z}^n:j_h \geq 0\
\mbox{for}\ 1 \leq h \leq n, \ j_1+j_2+\cdots+j_n=k\right\} \leq
(k+1)^{n-1}.
$$
 This completes the proof of
Lemma 4.2. $\endpf$

\medskip
To carry out the rapid iteration procedure, we need the following
estimates of the solution given by (\ref{2.10}) for the  equation
(\ref{4.6}).
\bigskip

{\bf Proposition 4.3:}\hspace{0.2cm}  {\it Suppose that
$w=|z|^2+E(z,\bar{z})$ is formally  equivalent to $M_\infty$ with
$E$ holomorphic over $\-{D_r}$ and  $Ord(E) \geq d$. Then the
solution given in (\ref{2.10})
satisfies the following estimates:\\
\begin{equation}
\begin{array}{cll}
|\hat f_h|_{\varrho},|\hat{g}|_{\varrho} &\leq&
 \frac{C(n) (2d)^{2n}\|E\|_r}{r-\varrho}(\frac{\varrho}{r})^{d-1}, \\
|\nabla \hat{f}_h|_{\varrho},|\nabla \hat{g}|_{\varrho} &\leq&
\frac{C(n)(2d)^{2n}\|E\|_r}{(r-\varrho)^3}(\frac{\varrho}{r})^{\frac{d-1}{2}},\\
|\hat{\varphi}|_{\varrho} &\leq&
\frac{(2d)^{2n}\|E\|_r}{(r-\varrho)^{2n}}(\frac{\varrho}{r})^{2d-2},
\label{4.12}
\end{array}
\end{equation}
where $C(n)=3^3n(n+1)2^{n+3}$. }\bigskip

{\it Proof of Proposition 4.3:} Notice that by the definition of
$\hat{f}$ given in (\ref{4.7}), we have $Ord(\hat{f})\ge d-1$ and
$deg_{wt}\leq 2d-4$. In terms of (\ref{2.10}), we can write, for
$2\leq h \leq n$,
$$
\hat{f}_h=A_1+A_2+A_3+A_4,
$$
where
\begin{equation*}
\begin{array}{l}
A_1=\sum\limits_{d\leq 2k+2 \leq
    2d-3}\frac{1}{2}z_h\left(-\sum_{j=2}^{h-1}Re(
    E_{(0,0)}^{(ke_1+e_j)}u^k)-2Re (
    E_{(0,0)}^{(ke_1+e_h)}u^k)\right),\\
A_2=\sum\limits_{i>h,d\leq 2k+2 \leq
2d-3}z_iE_{(e_i,e_h)}^{(ke_1)}u^k,\ A_3=\sum\limits_{|J|\geq 1,d
\leq |J|+2k+1 \leq
    2d-3}z^JE_{(J,e_h)}^{(ke_1)}u^k,\\
A_4=\sum\limits_{|I|\geq 1,d \leq |I|+2k \leq
    2d-3}z^{I+e_h}u^{k-1}\-{E^{(ke_1)}_{(0,I)}}+\sum\limits_{|I|\geq 1,d
    \leq |I|+2k+2 \leq
    2d-3}\sum\limits_{i=0}^{n-h-1}z^{I+e_h}u^{k}\-{E^{(ke_1+e_{n-i})}_{(0,I)}}\\
\hskip 1cm  -\sum\limits_{|I|\geq
    1,d \leq |I|+2k+2 \leq
    2d-3}z^{I+e_h}u^{k}\-{E^{(ke_1+e_h)}_{(0,I)}}\\
\hskip 0.4cm :=B_1+B_2+B_3.
\end{array}
\end{equation*}
By Lemma 4.2, we have, for $B_1$, the following\\
\begin{equation*}
\begin{array}{lll}
\|B_1\|_{\varrho}&=&\|\sum\limits_{|I|\geq 1,d \leq |I|+2k \leq
    2d-3}z^{I+e_h}u^{k-1}\-{E^{(ke_1)}_{(0,I)}}\|_{\varrho}\\
&\leq&\sum\limits_{|I|\geq 1,d \leq |I|+2k \leq
    2d-3}(R')^{I+e_h}(2{\varrho}^2)^{k-1}\frac{(k+2)^n\|E\|_r}{R^{I}\cdot(2r^2)^k}\\
&\le&\sum\limits_{|I|\geq 1,d \leq |I|+2k \leq
    2d-3}\left(\frac{\varrho}{r}\right)^{|I|+2k-1}\frac{(k+2)^n\|E\|_r}{2r}\\
&\leq&\sum\limits_{|I|\geq 1,d \leq |I|+2k \leq
    2d-3}\left(\frac{\varrho}{r}\right)^{|I|+2k-1}(2d)^n\|E\|_r\\
&\leq&\frac{(2d)^{2n}\|E\|_r}{r-\varrho}
\left(\frac{\varrho}{r}\right)^{d-1}.
\end{array}
\end{equation*}
Here and in what follows, we write
$R'=(2^{-\frac{n-2}{2}}\varrho,2^{-\frac{n-2}{2}}\varrho,2^{-\frac{n-3}{2}}\varrho\cdots,2^{-\frac{1}{2}}\varrho,\varrho)$.
Wee have also used the fact that
$$
\sharp\{(i_1,i_2,\cdots,i_n,k)\in \mathbb{Z}^{n+1}:i_h,k \geq 0 \
\mbox{for}\ 1\leq h \leq n,\sum_{h=1}^{n}i_h+2k=2d-1\}\leq (2d)^n.
$$
For $B_2$, we have
\begin{equation*}
\begin{array}{lll}
\|B_2\|_{\varrho}&=&\|\sum\limits_{|I|\geq 1,d
    \leq |I|+2k+2 \leq 2d-3}\sum\limits_{i=0}^{n-h-1}z^{I+e_h}u^{k}\-{E^{(ke_1+e_{n-i})}_{(0,I)}}\|_{\varrho}\\
&\leq&\sum\limits_{|I|\geq 1,d \leq |I|+2k+2 \leq
    2d-3}(R')^{I+e_h}(2{r'}^2)^{k}\cdot n
    \frac{2^n(k+2)^n\|E\|_r}{R^{I}\cdot(2r^2)^{k+1}}\\
&\le&\sum\limits_{|I|\geq 1,d \leq |I|+2k+2 \leq
    2d-3}\left(\frac{\varrho}{r}\right)^{|I|+2k+1}\frac{n2^n(k+2)^n\|E\|_r}{2r}\\
&\leq&\frac{n2^n(2d)^{2n}\|E\|_r}{r-\varrho}
\left(\frac{\varrho}{r}\right)^{d-1}.
\end{array}
\end{equation*}
Similarly, we have $\|B_3\|_{\varrho}\leq
\frac{2^n(2d)^{2n}\|E\|_r}{r-\varrho}
\left(\frac{\varrho}{r}\right)^{d-1}$. Thus we obtain:
$$
\|A_4\|_{\varrho}\leq \frac{n2^{n+1}(2d)^{2n}\|E\|_r}{r-\varrho}
\left(\frac{\varrho}{r}\right)^{d-1}.
$$
In the same manner, we have
$$
\left\|A_1\right\|_{\varrho},\left\|A_2\right\|_{\varrho},\left\|A_3\right\|_{\varrho}\leq\frac{n\cdot2^{n+1}(2d)^{2n}\|E\|_r}{r-\varrho}(\frac{\varrho}{r})^{d-1},
$$
Hence we get
$$
|\hat{f}_h|_{\varrho}\leq \frac{n \cdot
         2^{n+3}(2d)^{2n}\|E\|_r}{r-\varrho}(\frac{\varrho}{r})^{d-1}.
$$

Now letting $\tau=\frac{r+2\varrho}{3}$ and using the Cauchy estimates,
we have the following  estimate of the derivatives of $F$:
\begin{equation}
\begin{array}{l}
|(\hat{f}_h)'_{z_i}|_{\varrho} \leq \frac{
\tau|\hat{f}_h|_{\tau}}{(\tau-\varrho)^2}\leq
\frac{3^3n\cdot 2^{n+3}(2d)^{2n}\|E\|_r}{(r-\varrho)^3}(\frac{\varrho}{r})^{\frac{d-1}{2}},\\
|(\hat{f}_h)'_{w}|_{\varrho} \leq \frac{
2\tau^2|\hat{f}_h|_{\tau}}{(2\tau^2-2(\varrho)^2)^2}\leq
\frac{3^3n\cdot
2^{n+3}(2d)^{2n}\|E\|_r}{(r-\varrho)^3}(\frac{\varrho}{r})^{\frac{d-1}{2}}.
\label{2.4.13}
\end{array}
\end{equation}
Here we have used the fact that
\begin{equation}
(\frac{\tau}{r})^2\leq \frac{\varrho}{r}\ \mbox{for}\ \frac{1}{2}<
\varrho<\tau <r\leq 1\ ,\  \tau=\frac{r+2\varrho}{3}. \label{1.4.13}
\end{equation}
The inequality (\ref{2.4.13}) shows that $|\nabla
\hat{f}_h|_{\varrho}\leq \frac{3^3n(n+1)\cdot
2^{n+3}(2d)^{2n}\|E\|_r}{(r-r')^3}(\frac{\varrho}{r})^{\frac{d-1}{2}}$.
 The corresponding estimates on $\hat{f}_1$
and $\hat{g}$ can be achieved similarly.

We next  estimate  $\hat{\varphi}$. Notice that $-\hat{g}(z,u)
+2Re\left(\sum_{h=1}^{n}\bar{z_i}\hat{f}_i(z,u)\right)$ is only used
to cancel  terms of with weight $< 2d-2$ in $E$. By (\ref{4.6.0}),
we have the following:
\begin{equation*}
\begin{array}{l}
\| \hat{\varphi}\|_{\varrho}=\|\sum\limits_{t \geq 2d-2}E^{(t)}\|_{\varrho}\\
=\|\sum\limits_{|I|+|J| \geq
2d-2}a_{i_1\cdots i_nj_1\cdots j_n}z_1^{i_1}\cdots z_n^{i_n}{\overline{z_1}^{j_1}}\cdots{\overline{z_n}^{j_n}}\|_{\varrho}\\
\leq \sum\limits_{|I|+|J| \geq
2d-2}\|E\|_r(\frac{R'}{R})^{I+J}\\
\leq\sum\limits_{|I|+|J| = 2d-2,|K|,|L|\geq
0}\|E\|_r(\frac{\varrho}{r})^{|I|+|J|}\cdot
(\frac{\varrho}{r})^{k_1}\cdots(\frac{\varrho}{r})^{k_n}\cdot(\frac{\varrho}{r})^{l_1}\cdots(\frac{\varrho}{r})^{l_n}\\
\leq \sum\limits_{|I|+|J| = 2d-2}\|E\|_r
(\frac{\varrho}{r})^{2d-2}\cdot
(\frac{1}{1-\frac{\varrho}{r}})^{2n}\\
\leq
\frac{(2d)^{2n}\|E\|_r}{(r-\varrho)^{2n}}(\frac{\varrho}{r})^{2d-2}.
\end{array}
\end{equation*}
Here we have used the fact that
$$
\sharp \{(i_1,\cdots,i_n,j_1,\cdots,j_n)\in
\mathbb{Z}^{2n}:i_h,j_h\geq 0 \ \mbox{for}\ 1\leq h \leq
n,\sum_{h=1}^{n}(i_h+j_h)=k\} \leq (k+1)^{2n}.
$$
This finishes the proof of Proposition 4.3. $\endpf$\\

{\bf Proposition 4.4:}\hspace{0.2cm} {\it Let $E,r,\varrho, C(n)$ be
as in Proposition 4.3. Then there exists a constant $\delta>0$ such
that for
\begin{equation}
\frac{C(n)(2d)^{2n}\|E\|_r}{(r-\varrho)^3}(\frac{\varrho}{r})^{\frac{d-1}{2}}
<\delta, \label{4.13}
\end{equation}
$\Psi(z',w'):=\Theta^{-1}(z',w')$ is well defined in
$\-{\triangle_\sigma}$. Moreover, it  holds that
$\Psi(\triangle_{r'}) \subset \triangle_\sigma$,
$\Psi(\triangle_\sigma) \subset \triangle_\varrho$, $E'(z,\xi)$ is
holomorphic in $\-{\triangle_\sigma}$ and
\begin{equation}
\|E'\|_{r'} \leq C_d\|E\|_r^2+\widetilde{C_d}\|E\|_r. \label{4.14}
\end{equation}
Here
\begin{equation*}
\begin{array}{l}
C_d=\frac{(2n+1)\cdot3^3C(n)(2d)^{2n}}{(r-r')^3}(\frac{r'}{r})^{\frac{d-1}{4}}+
(\frac{r'}{r})^{d-1}n\cdot
\left(\frac{3C(n)(2d)^{2n}}{r-r'}\right)^2,\ \
\widetilde{C_d}=\frac{3^{2n}\cdot(2d)^{2n}}{(r-r')^{2n}}(\frac{r'}{r})^{d-1}.
\end{array}
\end{equation*}

 }
\medskip

{\it Proof of Proposition 4.4:} We need to show that for each
$(z',w') \in \-{\triangle_{\sigma}}$, we can uniquely solve the
system:
\begin{equation*}
\left\{
\begin{array}{l}
z'=z+\hat{f}(z,w) \\
w'=w+\hat{g}(z,w)
\end{array}
\right.
\end{equation*}
with $(z,w) \in \triangle_{\varrho}$. By (\ref{4.12}), choosing
$\delta$ sufficiently small  such that $|\nabla
\hat{f}|_{\varrho}+|\nabla \hat{g}|_{\varrho}<\frac{1}{2n+4}$ and
$|\hat{f}|_{\varrho}+|\hat{g}|_{\varrho}<\frac{1}{2n+4}.(r-\varrho)$ .
Define $(z^{[1]},w^{[1]})=(z',w') \in \triangle_{\sigma}$ and
$(z^{[j]},w^{[j]})$ inductively by
\begin{equation*}
\left\{
\begin{array}{l}
z^{[j+1]}=z'-\hat{f}(z^{[j]},w^{[j]}) \\
w^{[j+1]}=w'-\hat{g}(z^{[j]},w^{[j]}).
\end{array}
\right.
\end{equation*}

By a standard argument on the Picard iteration procedure, we
can get a unique $(z,w) \in \triangle_{\varrho}$ satisfying
$\Psi^{-1}(z,w)=(z',w')$, which gives that $\Psi(\triangle_\sigma)
\subset \triangle_\varrho$. Similarly, we have $\Psi(\triangle_{r'})
\subset \triangle_\sigma$. Hence we conclude that $E'$ is
holomorphic in $\triangle_\sigma$. Moreover,
\begin{equation}
\|E'\|_{r'} \leq \|Q\|_\sigma
\end{equation}
where
\begin{equation}
\begin{array}{l}
Q=\left(\hat{g}(z,\Phi)-\hat{g}(z,u)\right)-2Re\left(\sum\limits_{i=1}^{n}\overline{z_i}\Large(\hat{f}_i(z,\Phi)-\hat{f}_i(z,u)\Large)\right)
  -|\hat{f}(z,\Phi)|^2
  +\hat{\varphi}(z,\bar{z}).
\end{array}
\label{4.16}
\end{equation}
Notice that
\begin{equation}
\begin{array}{lll}
 |(\hat{g}(z,\Phi)-\hat{g}(z,u)|_{\sigma}
&\leq& |\nabla
  \hat{g}|_{\varrho} \cdot \|E\|_r \leq \frac{C(n)(2d)^{2n}\|E\|_r^2}{(r-\varrho)^3}(\frac{\varrho}{r})^{\frac{d-1}{2}}\\
&\leq&
\frac{3^3C(n)(2d)^{2n}\|E\|_r^2}{(r-r')^3}(\frac{r'}{r})^{\frac{d-1}{4}}.
\end{array}
\label{4.17}
\end{equation}
Here we have used the fact that $(\frac{\varrho}{r})^{2} <
\frac{r'}{r}$. (This can be achieved by the same token as for
(\ref{1.4.13}).)  We also have
\begin{equation}
\begin{array}{l}
|(\hat{f}_i(z,\Phi)-\hat{f}_i(z,u)|_{\sigma} \leq
\frac{3^3C(n)(2d)^{2n}\|E\|_r^2}{(r-r')^3}(\frac{r'}{r})^{\frac{d-1}{4}}\ \rm{for} \ 1\leq i \leq n.\\
|\hat{f}(z,\Phi)|^2_{\sigma} \leq n\cdot \left(\frac{C(n)
{(2d)}^{2n}\|E\|_r}{r-\sigma}(\frac{\sigma}{r})^{d-1}\right)^2 \leq
n \cdot \left(\frac{3C(n)
{(2d)}^{2n}\|E\|_r}{r-r'}\right)^2(\frac{r'}{r})^{d-1}\\
\|\hat{\varphi}\|_{\sigma}\leq
\frac{(2d)^{2n}\|E\|_r}{(r-\sigma)^{2n}}(\frac{\sigma}{r})^{2d-2}
\leq \frac{3^{2n}(2d)^{2n}\|E\|_r}{(r-r')^{2n}}(\frac{r'}{r})^{d-1}.
\end{array}\label{4.18}
\end{equation}
By (\ref{4.16})-(\ref{4.18}), we obtain:
\begin{equation*}
\begin{array}{l}
\|E'\|_{r'} \leq
\left\{\frac{(2n+1)\cdot3^3C(n)(2d)^{2n}}{(r-r')^3}(\frac{r'}{r})^{\frac{d-1}{4}}+
(\frac{r'}{r})^{d-1}n\cdot
\left(\frac{3C(n)(2d)^{2n}}{r-r'}\right)^2\right\}\|E\|_r^2+\frac{3^{2n}\cdot(2d)^{2n}}{(r-r')^{2n}}(\frac{r'}{r})^{d-1}\|E\|_r
\end{array}
\end{equation*}
This completes the proof of Proposition 4.4. $\endpf$
\bigskip

Now we turn to the proof of  Theorem 1. Set
$r_{\upsilon},\varrho_{\upsilon},\sigma_{\upsilon}$ as follows:
$$
r_{\upsilon}=\frac{1}{2}\left(1+\frac{1}{\upsilon+1}\right)\ ,\
\varrho_{\upsilon}=\frac{1}{3}(2r_{\upsilon}+r_{\upsilon+1})\ ,\
\sigma_{\upsilon}=\frac{1}{3}(2r_{\upsilon}+\varrho_{\upsilon}).
$$
We will apply the previous estimates with
$r=r_\upsilon,\varrho=\varrho_\upsilon,\sigma=\sigma_\upsilon,
r'=r_{\upsilon+1}, \Psi=\Psi_v, \cdots,$ with $v=0,1,\cdots,.$ Then
we have the following (see [(4.5), Moser]):

\begin{equation}
(r_\upsilon-r_{\upsilon+1})^{-1}=2( \upsilon+1)(\upsilon+2),\
\frac{r_{\upsilon+1}}{r_\upsilon}=1-\frac{1}{(\upsilon+2)^2} \label
{new-01}
\end{equation}

 Define a sequence of real
analytic submanifolds
$$
M_k\ :\ w=|z|^2+E_k(z,\bar{z})
$$
by $M_0=M$, $M_{\upsilon+1}=\Psi^{-1}_{\upsilon}(M_{\upsilon})$ for
all $\upsilon=0,1,2,\cdots$, where $\Psi_\upsilon$ is the
biholomorphic mapping taking $\triangle_{\sigma_{\upsilon}}$ into $
\triangle_{\varrho_{\upsilon}}$. And let
$$
d_{\upsilon}=Ord(E_{\upsilon})\ ,\ \Phi_{\upsilon}=\Psi_0 \circ
\Psi_1 \circ \cdots \circ \Psi_{\upsilon}.
$$
Since $s=\infty$, we find that
$$
Ord(E_{\upsilon})=d_{\upsilon}\geq 2^{\upsilon}+2\ \mbox{for}\
\upsilon \geq 0.
$$

We next state the following elementary fact:

\medskip {\bf Lemma 4.5}: Suppose that there is a constant C and number $a>1$ such that $d_v\ge Ca^v$.
Then for any integer $m_1,m_2,m_3>0$,
$$ \lim_{v\ra \infty}
v^{m_3}d_v^{m_1}(1-\frac{1}{v^{m_2}})^{d_v}=0.$$

\medskip

Then one can prove, by using (\ref{new-01}) and Lemma 4.5, that
$$
\lim_{\upsilon \rightarrow \infty}C_{d_{\upsilon}}=0\ ,\
\lim_{\upsilon \rightarrow \infty}\widetilde{C_{d_{\upsilon}}}=0.
$$
Hence $C_{d_{\upsilon}}$ and $\widetilde{C_{d_{\upsilon}}}$ are
bounded. Set $C_{d_{\upsilon}},\widetilde{C_{d_{\upsilon}}} <C$,
where $C$ is a fixed positive constant.
Also, one can verify that the hypothesis in (\ref{4.13}) holds for
all $\upsilon \geq 0$,
 by choosing $\eta^*_0=\|E_0\|_{r_0}$ sufficiently small.
 Indeed, we can even have
$\|E_\upsilon\|_{r_\upsilon}\leq \epsilon 2^{-\upsilon}$ for all
$\upsilon \geq 0$ and any given $1>\epsilon>0$.

Choose $N$ large enough such that $C_{d_{\upsilon}},
\widetilde{C_{d_{\upsilon}}}\leq \frac{1}{4}$ when $\upsilon \geq
N$. Suppose $C>1$ and choose $E_{0}$ such that
$\eta^*_0=\epsilon (2C)^{-2N}<1$. Then we have the following\\

(I) When $\upsilon \leq N$, we have
\begin{equation*}
\begin{array}{l}
\|E_\upsilon\|_{r_\upsilon} \leq
C(\|E_{\upsilon-1}\|_{r_{\upsilon-1}}+1)\|E_{\upsilon-1}\|_{r_{\upsilon-1}}
\leq 2C \cdot \|E_{\upsilon-1}\|_{r_{\upsilon-1}} \leq
(2C)^{\upsilon}\|E_0\|_{r_{0}} \leq \epsilon (2C)^{\upsilon-2N} \leq
\epsilon 2^{-N}.
\end{array}
\end{equation*}

(II) When $\upsilon > N$, we have
\begin{equation*}
\begin{array}{l}
\|E_\upsilon\|_{r_\upsilon} \leq \frac{1}{4}\cdot 2\cdot
\|E_{\upsilon-1}\|_{r_{\upsilon-1}} \leq
(\frac{1}{2})^{\upsilon-N}\|E_N\|_{r_N} \leq \epsilon 2^{-\upsilon}.
\end{array}
\end{equation*}
Now, choose $\epsilon$ sufficiently small. Then it follows from
(\ref{4.12}) and Proposition 4.4 that
$\|d\Psi_{\upsilon}^{-1}\|_{\triangle_{\varrho_{\upsilon}}}\leq
1+C_0 \|E_\upsilon\|_{r_\upsilon}\leq 1+C_0\epsilon2^{-\upsilon}$
for some constant $C_0$. Notice that $\Psi_\upsilon$ maps
$\triangle_{\sigma_\upsilon}$ into $\triangle_{\varrho_\upsilon}$.
By Cramer's rule, we have
$\|d\Psi_{\upsilon}\|_{\triangle_{\sigma_{\upsilon}}}\leq 1+\epsilon
C_12^{-\upsilon}$ for some constant $C_1$. Now the convergence of
$\Phi_\upsilon$ in $\triangle_{\frac{1}{2}}$ follows from the fact
that
$$
\Pi_{\upsilon=0}^{\infty}\|d\Psi_{\upsilon}\|_{\triangle_{\sigma_{\upsilon}}}
\leq \Pi_{\upsilon=0}^{\infty}(1+\epsilon C_12^{-\upsilon})<\infty,
$$
which completes the proof of Theorem 1. $\endpf$

\bigskip
{\bf Remark 4.6}: We notice that the formal map in Theorem 1 sending
$(M,0)$ to its quadric $(M_\infty, 0)$ may not be convergent
 as $aut_0(M_\infty)$ contains
many non-convergent elements. This is quite different from the
setting for CR manifolds, where formal maps are always convergent
under certain not too degenerate assumptions. We refer the reader to
the survey article [BER1] for  discussions and references on this
matter.

\bigskip
\noindent Xiaojun Huang (huangx@math.rutgers.edu),
Department of Mathematics, Rutgers University at New Brunswick, NJ
08903, USA;

\noindent
 Wanke Yin, School of Mathematical Sciences, Wuhan University, Wuhan 430072, P. R. China.
\end{document}